%% file: affine-rigidity-LC.tex
\documentclass[11pt,twoside,reqno,openany]{amsart}

\pdfoutput=1

\usepackage{amsbsy,amscd,amsfonts,amsmath,amssymb,amsthm,color,
fancybox,fancyhdr,footmisc,graphics,graphicx,ifthen,mathrsfs,
multicol,pdfpages,rotating,times,wasysym}
\usepackage[dvipsnames,svgnames,x11names]{xcolor}

\usepackage{pifont}

\usepackage[all]{xy}
\usepackage[utf8]{inputenc}
\usepackage[T1]{fontenc}
\usepackage{mathtools}
\newtagform{EngelLie}[\scriptstyle]{$}{$}
\makeatletter\newcommand{\leqnomode}{\tagsleft@true}
\newcommand{\reqnomode}{\tagsleft@false}\makeatother

\input macros.tex


\begin{document}

\setcounter{section}{0}

$\:$


\begin{center}

{\large\bf On Degenerate Affine 
and Cauchy-Riemann Geometries\footnotemark[1]}


\label{affine-rigidity-without-integration}

\bigskip\bigskip

Jo\"el~{\sc Merker}\footnotemark[2]

\footnotetext[1]{\,
This work was supported
in part by the Polish National Science Centre (NCN) 
via the grant number 2018/29/B/ST1/02583.}

\footnotetext[2]{\,
Laboratoire de Math\'ematiques d'Orsay, CNRS, 
Universit\'e Paris-Saclay, 91405 Orsay Cedex, France.
{\bf joel.merker@universite-paris-saclay.fr}}

\end{center}\bigskip

\begin{center}
\begin{minipage}[t]{12.25cm}
\parindent 0.53cm
\scriptsize
\noindent
{\sc Abstract}.
Tube (real) hypersurfaces $M = H \times i\, \R^\NN$ in the complex
space $\C^\NN$, with $H \subset \R^\NN$ an $(\NN-1)$-dimensional
hypersurface, have dimension $2\NN-1$, but depend {\em in fine} only
on $\NN-1$ variables, since they are invariant under the imaginary
translations $\zaux \longmapsto \zaux + i\, c$, with $c \in
\R^\NN$. Their CR geometry is easier to understand, as for instance in
$\C^2$, Cartan's primary invariant $\Iaux_{\sf Cartan}$ has $5$ terms
when $M$ is tube, but $> {\bf 10^6}$ terms in the general case
(Merker-Sabzevari, Abel Symposium 2015).

Tubes $M = H \times i\, \R^\NN$ can be classified under the two
(local) Lie groups:
\[
\aligned
\Affsmall^\tubesmall(\C^\NN)
&
\,:=\,
\big\{
\zaux\longmapsto A\,\zaux+b+i\,c
\colon\,\,
A\in\GLsmall(\R^\NN),\,\,
b\in\R^\NN,\,\,
c\in\R^\NN
\big\},
\\
\Biholsmall(\C^\NN)
&
\,:=\,
\big\{
\zaux\longmapsto\zaux'(\zaux)
\colon\,
\text{locally biholomorphic}
\big\},
\endaligned
\]
of respective dimensions $\NN^2 + 2\NN < \infty$.  Then the (partial)
classifications known for Levi nondegenerate $M$, {\em i.e.}  for $H$
having nondegenerate Hessian, differ significantly, whatever
the signature is.

Dadok-Yang in 1985, and then Isaev in a Springer monograph of 2011,
considered {\sl spherical tubes} in $\C^{n+1} \ni (z_1, \dots, z_n,
w)$, namely tubes biholomorphic to $\Resmall\, w = \pm (\Resmall\,
z_1)^2 \pm \cdots \pm (\Resmall\, z_n)^2$, determined equivalence
classes under $\Affsmall^\tubesmall(\C^\NN)$, and obtained complete
results in signatures $(p, n-p)$ with $p = 0, 1, 2$ [in $\C^2$,
sphericity is equivalent to $\Iaux_{\sf Cartan} \equiv 0$].

For Levi {\em degenerate} (non Levi flat) tubes, similar
classification ramifications have been expected, and the recently much
studied class $\mathfrak{C}_{2,1}$ of $2$-nondegenerate constant Levi
rank $1$ hypersurfaces $M^5 \subset \C^3 \ni (z, \zeta, w)$ has been
considered by Isaev as a {\em test-case}. In this context, `{\sl
sphericity}', namely biholomorphic equivalence to the
Gaussier-Merker (maximally homogeneous) model, which can be graphed in
tube form as $\Resmall\, w = (\Resmall\, z)^2 \big/ (1 - \Resmall\,
\zeta)$, was characterized by Pocchiola as equivalent to the vanishing
of two invariants $0 \equiv \Waux_0 \equiv \Jaux_0$.

In {\em J. Differential Geom.} {\bf 104} (2016), 111--141, inspired by
Pocchiola, Isaev showed as a main\,\,---\,\,{\em
unexpected!}\,\,---\,\,theorem that any $\mathfrak{C}_{2,1}$ tube
hypersurface $M^5 = S^2 \times i\, \R^3$ which is biholomorphic to
$\Resmall\, w = (\Resmall\, z)^2 \big/ (1 - \Resmall\, \zeta)$ is in
fact already {\em affinely} equivalent to it. So, no classification
ramification occurs.

We provide a quicker proof, of length $< 2$ pages, which uses only
(straightforward) differential algebra\,\,---\,\,no integration.  We
also discuss smoothness: $\mathcal{C}^\omega$, $\mathcal{C}^\infty$,
$\mathcal{C}^5$.  Beyond, we explore the tight relationships between
Affine geometry and Cauchy-Riemann geometry, and we raise a few
accessible problems.
\end{minipage}
\end{center}

\medskip

\hfill
{\footnotesize\sf\em
Dedicated to the memory of Alexander 
Isaev\textsuperscript{$\dag$}}

\Section{\bf Introduction}
\label{introduction-affine-rigidity}
\HEAD{{\ref{introduction-affine-rigidity}}.~{\sf
Introduction}
}{
Jo\"el {\sc Merker}, 
Paris-Saclay University, Orsay, France}

The study of {\sl tube domains} $\Omega = D \times i\, \R^\NN$ in
$\C^\NN$, with $D \subset \R^\NN$ nonempty connected open set, $\NN
\geqslant 2$, is a classical subject in {\sl Several Complex
Variables}, which goes back (at least) to the beginning of the
20\textsuperscript{th} century.  This is exemplified by the celebrated
{\sl Bochner Tube Theorem} ({\cite[Thm.~2.5.10]{Hormander-1966}}),
which states that $\Omega$ has the simple envelope of holomorphy
$\widetilde{\Omega} := D_{\sf cvx} \times i\,\R^\NN$ , where $D_{\sf
cvx}$ is the {\em real} convex hull of $D$, {\em see}
also~{\cite{Merker-Porten-2006}} and~{\cite[IV.~12]{Merker-2017}}.
Naturally, {\sl tube hypersurfaces} $M^{2\NN-1} = H^{\NN-1} \times i\,
\R^\NN$ arise as boundaries of tube domains $D^\NN \times i\, \R^\NN$,
with the {\sl base} real hypersurface $H^{\NN-1} \subset \R^\NN$ being
the boundary of $D^\NN$.

Not only envelopes of holomorphy, but also differential invariants of
Cauchy-Riemann (CR) real hypersurfaces $M^{2\NN-1} \subset \C^\NN$,
and their Lie symmetry groups as well, become more
tractable\big/visible when a tube structure is supposed.  For
instance, the familiar unit ball $\big\{ \vert \zaux_1 \vert^2 +
\cdots + \vert \zaux_\NN \vert^2 < 1 \big\}$ in $\C^\NN \ni (\zaux_1,
\dots, \zaux_\NN)$ is biholomorphically equivalent to $D^\NN \times
i\, \R^\NN$ where:
\leqnomode\usetagform{default}
\begin{align}
\label{unbounded-sphere}
D^\NN
\colon\ \ \ \ \
\Re\,\zaux_\NN
\,>\,
(\Re\,\zaux_1)^2
+\cdots+
(\Re\,\zaux_{\NN-1})^2,
\end{align}
hence $D^\NN \times i\, \R^\NN$ is invariant through translations
along the {\em imaginary} axes. 
More generally,
Siegel found it often convenient to realize many symmetric
domains as tubes, {\em cf.}~{\cite{Satake-1980}},
and {\em see} also~{\cite{Eastwood-Ezhov-2004,
Kruzhilin-Soldatkin-2004, Kruzhilin-Soldatkin-2006}}
for classifications of tube {\em domains}.

Thus, the property that makes tube CR hypersurfaces $M^{2\NN-1} =
H^{\NN-1} \times i\, \R^\NN \subset \C^\NN$ interesting from the
complex-geometric point of view is that they possess an
$\NN$-dimensional commutative group of holomorphic symmetries, namely
the group of translations $\zaux \longmapsto \zaux + i\, c$, 
with arbitrary $c \in \R^\NN$.  

To be specific, introduce the (local) Lie groups:
\[
\aligned
\Aff(\R^\NN)
&
\,:=\,
\big\{
\xaux\longmapsto A\,\xaux+b
\colon\,\,
A\in\GL(\R^\NN),\,\,
b\in\R^\NN
\big\},
\\
\Aff^\tubesmall(\C^\NN)
&
\,:=\,
\big\{
\zaux\longmapsto A\,\zaux+b+i\,c
\colon\,\,
A\in\GL(\R^\NN),\,\,
b\in\R^\NN,\,\,
c\in\R^\NN
\big\},
\\
\Aff(\C^\NN)
&
\,:=\,
\big\{
\zaux\longmapsto\mathcal{A}\,\zaux+\beta
\colon\,\,
\mathcal{A}\in\GL(\C^\NN),\,\,
\beta\in\C^\NN
\big\},
\\
\Bihol(\C^\NN)
&
\,:=\,
\big\{
\zaux\longmapsto\zaux'(\zaux)
\colon\,
\text{locally biholomorphic}
\big\},
\endaligned
\]
which satisfy:
\[
\Aff(\R^\NN)
\,\hookrightarrow\,
\Aff^\tubesmall(\C^\NN)
\,\subset\,
\Aff(\C^\NN)
\,\subset\,
\Bihol(\C^\NN),
\]
Their {\em real} dimensions are $\NN^2 + \NN <
\NN^2 + 2\NN
< (2\NN)^2 + 2\NN < \infty$. 
For a tube $M = H \times i\, \R^\NN$,
introduce also:
\[
\aligned
\Sym^\Affsmall(H)
&
\,:=\,
\big\{
g\in\Aff(\R^\NN)
\colon\,
g(H)
\subset
H
\big\},
\\
\Sym^{\Affsmall^{\tubesmall}}
(M)
\,=\,
\Sym^{\Affsmall^{\tubesmall}}
\big(H\times i\,\R^\NN\big)
&
\,:=\,
\big\{
h\in\Aff^\tubesmall(\C^\NN)
\colon\,
h(M)
\subset
M
\big\},
\\
\Sym^{\Biholsmall}
(M)
\,=\,
\Sym^{\Biholsmall}
\big(H\times i\,\R^\NN\big)
&
\,:=\,
\big\{
h\in\Bihol(\C^\NN)
\colon\,
h(M)
\subset
M
\big\}.
\endaligned
\]
the inclusion symbol `$\subset$' being understood
in a local sense, and the group
elements $g$, $h$ being close to the identity.

Thus, one always has:
\[
i\,\R^\NN
\,=\,
\big\{
\zaux
\longmapsto
\zaux+i\,c
\big\}
\,\subset\,
\Sym^{\Affsmall^{\tubesmall}}
\big(H\times i\,\R^\NN\big),
\]
and in fact, one can convince oneself that:
\[
\Sym^{\Affsmall^{\tubesmall}}
\big(H\times i\,\R^\NN\big)
\,=\,
\Sym^\Affsmall(H)
\ltimes
i\,\R^\NN.
\]
But the determination of $\Sym^\Affsmall(H)$ for any hypersurface
$H^{\NN-1} \subset \R^\NN$ of any dimension is a problem of infinite
complexity, like that of classifying all abstract Lie algebras.

At least, it is clear that:
\leqnomode\usetagform{default}
\begin{align}
\label{Aff-included-Bihol}
\Sym^\Affsmall(H)
\,\hookrightarrow\,
\Sym^{\Affsmall^{\tubesmall}}
(M)
\,\subset\,
\Sym^{\Biholsmall}
(M),
\end{align}
the last inclusion being in general {\em strict}, 
because of the jump to $\infty$-dimension.
In the present paper we look at hypersurfaces from both the
affine-geometric and CR-geometric points of view.

\begin{Problem}
\label{Pbm-Aff-Bihol-equivalences}
{\sl 
Study and classify real hypersurfaces
$H \subset \R^\NN$ under
the {\em finite-dimensional} Lie group $\Aff(\R^\NN)$,
for instance in `small' accessible 
dimensions\footnote{\,
For the fields 
$\K = \C$ or $\R$, the elementary classification of affinely
homogeneous curves $C^1 \subset \K^2$ appears
in~{\cite{Schirokow-Schirokow-1959}}.  The classification of surfaces
$S^2 \subset \K^3$ was terminated
in~{\cite{Doubrov-Komrakov-Rabinovich-1996,
Abdalla-Dillen-Vrancken-1997}}, and reobtained
in~{\cite{Eastwood-Ezhov-1999}}.  A classification of affinely
homogeneous surfaces in $\R^4$ and in $\P^4(\R)$ having symmetry
algebras of dimension $\geqslant 4$ appears
in~{\cite[p.~37]{Doubrov-Komrakov-1998}}.  We ignore whether $3$-folds
$V^3 \subset \K^4$ have been affinely classified, even partly.}.
Then consider the associated tubes $H \times i\, \R^\NN \subset
\C^\NN$ and determine how affine equivalence classes {\rm merge}
(become equivalent) under the {\em infinite-dimensional} local Lie
group $\Bihol(\C^\NN)$.
}
\end{Problem}

The {\sl reverse} problem starts by classifying 
tubes $H \times i\, \R^\NN$ modulo $\Bihol (\C^\NN)$, 
hopefully getting a list,
before attempting to split further each obtained class as
several classes under the smaller group $\Aff (\R^\NN)$.

Indeed, because of strict inclusion
in~{\eqref{Aff-included-Bihol}}, affine classification lists
should in general contain 
{\em more items} that biholomorphic classification
lists.  For instance, in $\C^2 \ni (z, w) = (x+iy, u+iv)$, it is known
({\cite{Loboda-1996, Merker-2015}}) that a tube $\{ u = F(x)\}$ with
$F_{xx} \neq 0$ is holomorphically equivalent to the sphere $\{ \Re\,
w = (\Re\, z)^2 \}$ if and only if $F$ satisfies a certain
6\textsuperscript{th} order ordinary differential equation
shown in Corollary~{\ref{Cor-spherical-tube-C2}}.
One verifies that $\{ u = e^x \}$
satisfies this {\sc ode}, whereas it clearly is {\em not} affinely
equivalent to $\{ u = x^2 \}$.

Given a hypersurface $H = \{ \rho(\xaux_1, \dots, \xaux_\NN) = 0
\big\}$ with $d\rho \neq 0$ on $\{ \rho = 0 \}$, its Hessian
determinant:
\[
\Hessian(\rho)
\,:=\,
\left\vert\!
\begin{array}{cccc}
0 & \rho_{\xaux_1} & \cdots & \rho_{\xaux_\NN}
\\
\rho_{\xaux_1} & \rho_{\xaux_1\xaux_1} & \cdots &
\rho_{\xaux_1\xaux_\NN}
\\
\vdots & \vdots & \ddots & \vdots
\\
\rho_{\xaux_\NN} & \rho_{\xaux_\NN\xaux_1} & \cdots &
\rho_{\xaux_\NN\xaux_\NN}
\end{array}
\!\right\vert
\]
is a relative invariant under affine transformations, and one calls
$H$ {\sl nondegenerate} when $\Hessian(\rho) (\xaux) \neq 0$ at every 
point $\xaux \in H$.

Similarly, using the operators $\partial_{\zaux_k} := \frac{1}{2}\,
\big( \partial_{\xaux_k} - i\, \partial_{\yaux_k} \big)$ and
$\partial_{\overline{\zaux}_k} := \frac{1}{2}\, \big(
\partial_{\xaux_k} + i\, \partial_{\yaux_k} \big)$, a real
hypersurface $M \subset \C^\NN$ defined implicitly as $\big\{
\rho(\xaux_1, \yaux_1, \dots, \xaux_\NN, \yaux_\NN) = 0 \big\}$ is
{\sl Levi nondegenerate} 
{\cite{Merker-Porten-2006}}
if:
\[
0
\,\neq\,
\Levi(\rho)
\,:=\,
\left\vert\!
\begin{array}{cccc}
0 & \rho_{\zaux_1} & \cdots & \rho_{\zaux_\NN}
\\
\rho_{\overline{\zaux}_1} & 
\rho_{\overline{\zaux}_1\zaux_1} & \cdots &
\rho_{\overline{\zaux}_1\zaux_\NN}
\\
\vdots & \vdots & \ddots & \vdots
\\
\rho_{\overline{\zaux}_\NN} & 
\rho_{\overline{\zaux}_\NN\zaux_1} & \cdots &
\rho_{\overline{\zaux}_\NN\zaux_\NN}
\end{array}
\!\right\vert.
\]
Again, $\Levi(\rho)$ is a relative invariant under
$\Bihol(\C^\NN)$.
Since $\rho_{\zaux_k} = \frac{1}{2}\, \rho_{\xaux_k}$ when $\rho$
depends only on $\xaux_1, \dots, \xaux_\NN$ (and not on $\yaux_1,
\dots, \yaux_\NN$), a tube $M = H \times i\, \R^\NN$ is Levi
nondegenerate if and only if its base $H$ is nondegenerate.  Most
publications in Affine geometry and in CR geometry were done under
such nondegeneracy assumptions, 
further taking account of the {\em signatures}
of $\Hessian (\rho)$ and of $\Levi (\rho)$.

Setting $\NN =: n + 1$, and working locally in $\R^{n+1} \ni (x_1,
\dots, x_n, u)$, given a graphed hypersurface $H = \big\{ u = F(x_1,
\dots, x_n) \big\}$, the Hessian can better be viewed as the $n \times
n$ matrix $\big( F_{x_j x_k} \big)$.  After some affine transformation
centered at some point of $H$ which becomes the origin, one can make
$u = x_1^2 + \cdots + x_p^2 -
x_{p+1}^2 - \cdots - x_n^2 + {\rm O}(3)$, so
that the Hessian signature at the
origin (and nearby) reads as $(p, n-p)$ for a certain integer $0
\leqslant p \leqslant n$.

Problem~{\ref{Pbm-Aff-Bihol-equivalences}}
being probably too wide, even in `small' dimensions, it has been
`restricted' by Dadok-Yang~{\cite{Dadok-Yang-1985}} to the class of
{\sl spherical} CR hypersurfaces $M^{2\NN - 1} \subset \C^\NN$, namely
those that are locally biholomorphic to the `sphere' $\big\{ \Re\, w =
(\Re\, z_1)^2 + \cdots + (\Re\, z_n)^2 \big\}$\,\,---\,\,in 
its unbounded
representation~{\eqref{unbounded-sphere}}.  Thus, one considers a
unique class under $\Bihol (\C^{n+1})$, the most CR-symmetric one, and
the {\em reverse} Problem~{\ref{Pbm-Aff-Bihol-equivalences}} is to
determine {\em all} equivalences classes under the {\em smaller} group
$\Aff^{\tubesmall} (\C^{n+1})$.

Dadok-Yang~{\cite{Dadok-Yang-1985}} were able to settle this
sub-problem in any CR dimension $n \geqslant 1$, and they showed that
all $M = H \times i\,\R^{n+1}$, with $H$ having nondegenerate {\em
positive} Hessian, {\em i.e.} of signature $(n, 0)$, are in the
mutually exclusive classification list reformulated
in~{\cite[p.~93]{Isaev-2011}}:

\smallskip\noindent{\bf (a)}\,
$\big\{ u = e^{x_1} + \cdots + e^{x_\nu} + 
x_{\nu+1}^2 + \cdots + x_n^2 \big\}$ for any
$0 \leqslant \nu \leqslant n$;

\smallskip\noindent{\bf (b)}\,
$\big\{ u = \arcsin\, (e^{x_1} + \cdots + e^{x_n}) \big\}$;

\smallskip\noindent{\bf (c)}\,
$\big\{ u = \log\, (1 - e^{x_1} - \cdots - e^{x_n}) \big\}$.

\smallskip

Turning to other signatures, Isaev attempted to fully classify
spherical $M^{2n+1} \subset \C^{n+1}$ under $\Aff^{\tubesmall}
(\C^{n+1})$.  Devoting a whole monograph~{\cite{Isaev-2011}} to this
(unexpectedly wide) sub-problem, he obtained complete results for $p =
0, 1, 2$ and partially classified 
collections for any signature $(p, n-p)$.

In all cases encountered, each equivalence class under
$\Bihol(\C^\NN)$ did {\sl split up} into {\em several} inequivalent
classes under $\Aff^{\tubesmall} (\C^\NN)$.
This raised an intriguing

\begin{Question}
{\sl Given two tubes $M = H \times i\,\R^\NN$ and
$M' = H' \times i\, {\R'}^\NN$, 
can it happen that:}
\[
M
\overset{\Biholsmall}{\,\,\cong\,\,}
M'
\ \ \ \ \ \ \ \ \ \ \ \ \ \ \ \ \ \ \ \
\overset{\text{\red{\bf ?}}}{\Longrightarrow}
\ \ \ \ \ \ \ \ \ \ \ \ \ \ \ \ \ \ \ \
M
\overset{\Affsmall}{\,\,\cong\,\,}
M'
\ \ \ \ \ \ \ \ \ \ \
\text{\sl or equivalently}
\ \ \ \ \ \ \ \ \ \ \ 
H
\overset{\Affsmall}{\,\,\cong\,\,}
H'
\text{\red{\bf ?}}
\]
\end{Question}

Quite unexpectedly, a positive answer to this question was discovered,
but in the context of {\em Levi degenerate} CR hypersurfaces. 
Indeed,
the following rigidity result, obtained by Isaev in~{\cite{Isaev-2016}}
is in stark contrast to the Levi nondegenerate case, where
the CR-geometric and affine-geometric classifications significantly
differ.
We explain undefined terms below.

\begin{Theorem}
\label{Thm-main}
For a $\mathcal{C}^\omega$ tube CR hypersurface $M^5 = S^2 \times i\,
\R^3$ in $\C^3$ which is everywhere $2$-nondegenerate and of constant
Levi rank $1$, the following two conditions are equivalent.

\smallskip{\bf (i)}\,
Its base surface $S^2$ is equivalent within $\R^3$, under the
$12$-dimensional group $\Aff_3(\R)$, to $\big\{ u = \frac{x^2}{1-y}
\big\}$.

\smallskip{\bf (ii)}\,
$M^5$ itself is equivalent within $\C^3$, 
under the $\infty$-dimensional
group $\Bihol_3(\C)$, to $\big\{ \Resmall\, w = \frac{(\Resmall\,
z)^2}{1-\Resmall\, \zeta} \big\}$.

\end{Theorem} 

Of course, {\small\bf (i)} $\Longrightarrow$ {\small\bf (ii)},
hence the reverse implication is the main thing.
Denoting coordinates on $\C^3$ by $\big(z, \zeta, w = u + iv\big)$,
with $M$ graphed as:
\[
\Re\,w
\,=\,
F\big(\Re\,z,\,\Re\,\zeta\big),
\]
the assumptions of $2$-nondegeneracy and of constant
Levi rank $1$ read as follows, 
{\em cf.}~{\cite{Chen-Foo-Merker-Ta-2020, 
Foo-Merker-Ta-2020}},
with
$x := \Re\, z$ and $y := \Re\, \zeta$:
\[
0
\,\neq\,
\left\vert\!
\begin{array}{cc}
F_{xx} & F_{xy}
\\
F_{xxx} & F_{xxy}
\end{array}
\!\right\vert
\ \ \ \ \ \ \ \ \ \ \ \ \ \ \ \ \ \ \ \
\text{and}
\ \ \ \ \ \ \ \ \ \ \ \ \ \ \ \ \ \ \ \
F_{xx}
\,\neq\,
0
\,\equiv\,
\left\vert\!
\begin{array}{cc}
F_{xx} & F_{xy}
\\
F_{xy} & F_{yy}
\end{array}
\!\right\vert.
\]
In Theorem~{\ref{Thm-main}} above, $\big\{ u = \frac{x^2}{1-y}
\big\}$ is the tube representation of
(is biholomorphically equivalent to)
the {\em maximally
homogeneous} Gaussier-Merker model~{\cite{Gaussier-Merker-2004}}:
\[
\Re\,w
\,=\,
\big(
z\overline{z}+\tfrac{1}{2}\,z^2\overline{\zeta}
+\tfrac{1}{2}\,\overline{z}^2\zeta
\big)
\big/
\big(
1-\zeta\overline{\zeta}
\big),
\]
whose biholomorphic automorphisms group is $10$-dimensional,
isomorphic to $\SO_{3,2}(\R)$, {\em see}~{\cite[p.~7
sq.]{Chen-Foo-Merker-Ta-2020}} for Lie group considerations, not
useful here.

Isaev's original proof of Theorem~{\ref{Thm-main}}
in~{\cite{Isaev-2016}} is advanced, demanding, and requires to {\em
integrate} certain {\sc pde}s. In this paper, we propose a (much)
quicker proof, of length $< 2$ pages, which uses only
(straightforward) differential algebra\,\,---\,\,no integration.

Pocchiola's Ph.D. was inspirational to~{\cite{Isaev-2016}}, since
after the {\footnotesize\sf arxiv}
prepublication~{\cite{Pocchiola-2013}}, Isaev re-obtained CR
invariants which characterized local biholomorphic equivalence of {\em
tubes} $\big\{ \Re\, w = F(\Re\, z, \, \Re\, \zeta) \big\}$ to the
Gaussier-Merker model, or to $\big\{ \Re\,w = \frac{(\Resmall
z)^2}{1-\Resmall \zeta} \big\}$.  However, we would like to mention
that Pocchiola's characterization of equivalence to $\big\{ \Re\,w =
\frac{(\Resmall z)^2}{1-\Resmall \zeta} \big\}$ was valid for any
general graphed hypersurface (not necessarily tube):
\[
\Re\,w
\,=\,
F\big(
\Re\,z,\,\Im\,z,\,
\Re\,\zeta,\,\Im\,\zeta,\,
\Im\,w
\big),
\]
with a graphing function depending on $5$, instead of $2$,
real variables, which raises up (significantly)
the level of computational complexity.

Sections~{\ref{Pocchiola-invariants-C21-M5-C3}} 
and~{\ref{W-aff-J-aff}}
discusses these aspects, assuming that the reader
already got acquainted a bit with CR geometry. Then
Section~{\ref{affine-rigidity-elimination}} 
presents the short proof of
Theorem~{\ref{Thm-main}}.  
Section~{\ref{smoothness-assumption-improvement}} 
explains why this statement
holds true with $M$ of class $\mathcal{C}^5$ (or
$\mathcal{C}^\infty$) instead of $\mathcal{C}^\omega$, {\em
cf.}~{\cite{Isaev-Merker-2019}}.


\medskip\noindent{\bf Acknowledgments.}
In February 2019, Alexander Isaev 
visited Orsay University, gave impetus, and fostered
with breadth exciting exchanges about relationships between CR
geometry and Affine geometry. 
For his generosity,
let him be thanked with
sheer gratitude,
from eternal ether. 

The inspiring relationships between Affine geometry
and CR geometry were then explored in further 
works~{\cite{Foo-Merker-2019,
Foo-Merker-Ta-2019,
Chen-Merker-2020,
Chen-Foo-Merker-Ta-2020,
Foo-Merker-Ta-2020,
Merker-Nurowski-2019,
Merker-Nurowski-2020,
Merker-Nurowski-2021, Merker-Nurowski-2022}}, 
conducted under 
the guidance of Pawe{\l} Nurowski.
The-Anh Ta read carefully the manuscript.

Lastly, we are glad to `advertise' that 
Question~{\ref{Q-differential-invariants-surfaces}}, appearing
in the {\small\tt arxiv.org} prepublication of this article, was
recently solved by \"Orn Arnaldsson and Francis Valiquette
in~{\cite{Arnaldsson-Valiquette-2020}}.

\Section{\bf Prologue: The $\C^2$ case, a Sketch}
\label{prologue-C2-case}
\HEAD{{\ref{prologue-C2-case}}.~{\sf Prologue: 
The $\C^2$ case, a Sketch}
}{
Jo\"el {\sc Merker}, 
Paris-Saclay University, Orsay, France}

Without complete details, let us briefly explain
why biholomorphic equivalence and affine equivalence
of a tube $\{ u = F(x) \}$ in $\C^2$
to $\{u' = (x')^2\}$ {\em differ}.

Let a hypersurface
$M^3 \subset \C^2$ be given in coordinates 
$(z,w) = (x+i\,y,\,u+i\,v)$ as a real $\mathcal{C}^\omega$ graph:
\[
u
\,=\,
F\big(x,y,v\big).
\]
Assume that $M$ is Levi nondegenerate. 
In the intrinsic coordinates $(x, y, v)$ on $M$, two
generators of $T^{1,0}M$ and $T^{0,1}M$ are (detailed explanation
appears in~{\cite[2.1]{Merker-Pocchiola-2018}}):
\[
\mathcal{L}
\,:=\,
\frac{\partial}{\partial z}
+
A\,\frac{\partial}{\partial v}
\ \ \ \ \ \ \ \ \ \ \ \ \ \ \ \ \ \
\text{and}
\ \ \ \ \ \ \ \ \ \ \ \ \ \ \ \ \ \
\overline{\mathcal{L}}
\,:=\,
\frac{\partial}{\partial\overline{z}}
+
\overline{A}\,\frac{\partial}{\partial v},
\]
where:
\[
A
\,:=\,
-\,i\,
\frac{F_{z}}{1+i\,F_v}.
\]
The Levi nondegeneracy assumption is equivalent to the
everywhere nonvanishing of the {\sl Levi factor:}
\[
\laux
\,:=\,
i\,
\Big(
\overline{A}_z
+
A\,\overline{A}_v
-
A_{\overline{z}}
-
\overline{A}\,A_v
\Big)
\,\,\neq\,\,
0.
\]
Introduce also a function whose complete expansion in terms
of $J_{x,y,v}^3 F$ is one page long
({\em cf.} {\cite[p.~42]{Aghasi-Merker-Sabzevari-2011}}):
\[
\Paux
\,:=\,
\frac{\laux_z+A\,\laux_v-\laux\,A_v}{\laux}.
\]

\begin{Theorem}
\label{Thm-sphericity-M3-C2}
{\rm\cite{Merker-2015}}
A Levi nondegenerate $\mathcal{C}^\omega$ local hypersurface
$M^3 \subset \C^2$ is locally biholomorphically equivalent to
the tube representation of the unit sphere:
\[
M
\overset{\Biholsmall}{\,\,\cong\,\,}
\big\{
u
=
x^2
\big\},
\]
if and only if:
\reqnomode\usetagform{EngelLie}
\begin{align}
0
\,\equiv\,
\Iaux_{\Cartansmall}
&
\,:=\,
-\,2\,\overline{\mathcal{L}}
\big(\mathcal{L}\big(\overline{\mathcal{L}}
(\overline{\Paux})
\big)\big)
+
3\,\overline{\mathcal{L}}\big(
\overline{\mathcal{L}}\big(
\mathcal{L}(\overline{\Paux})\big)\big)
-
7\,
\overline{\Paux}\,
\overline{\mathcal{L}}\big(
\mathcal{L}(\overline{\Paux})\big)
\,+
\notag
\\
&
\ \ \ \ \ \ \
+
4\,\overline{\Paux}\,\mathcal{L}\big(
\overline{\mathcal{L}}(\overline{\Paux})\big)
-
\mathcal{L}\big(\overline{\Paux}\big)\,
\overline{\mathcal{L}}\big(\overline{\Paux}\big)
+
2\,\overline{\Paux}\,\overline{\Paux}\,
\mathcal{L}\big(\overline{\Paux}\big).
\tag{\qed}
\end{align}
\end{Theorem}

Unfortunately, the real and imaginary parts of $\Iaux_{\Cartansmall}$
contain $> {\bf 10^6}$ differential monomials in $J_{x,y,v}^6F$, {\em
cf.}~{\cite[p.~178]{Merker-2015}}.  But when $M = \big\{ u = F(x)
\big\}$ is tube, the $1$ page long expression of $\Paux$ contracts as:
\[
\Paux
\,=\,
\frac{1}{2}\,
\frac{F_{xxx}}{F_{xx}}
\,=\,
\overline{\Paux}.
\]
Since $\Paux$, and $\laux$ as well, are functions 
of only $x$, hence are independent of $v$, 
the $(1,0)$ and $(0,1)$ differentiation operators
$\mathcal{L}$ and $\overline{\mathcal{L}}$ act
on them simply as $\frac{1}{2}
\frac{\partial}{\partial x}$.
Then the formula of Theorem~{\ref{Thm-sphericity-M3-C2}} 
becomes expandable.

\begin{Corollary}
\label{Cor-spherical-tube-C2}
When the hypersurface $M^3 \subset \C^2$ is a tube graphed
as $\big\{ u = F(x)\big\}$, it holds:
\reqnomode\usetagform{EngelLie}
\begin{align}
\Iaux_{\Cartansmall}
&
\,=\,
\frac{1}{16}\,
\Big\{
\big(F_{xx}\big)^3\,
F_{xxxxxx}
-
7\,\big(F_{xx}\big)^2\,
F_{xxx}\,F_{xxxxx}
-
4\,\big(F_{xx}\big)^2\,
\big(F_{xxxx}\big)^2
\,+
\notag
\\
&
\ \ \ \ \ \ \ \ \ \ \ \ \ \
+
25\,F_{xx}\,\big(F_{xxx}\big)^2\,F_{xxxx}
-
15\,\big(F_{xxx}\big)^3
\Big\}.
\tag{\qed}
\end{align}
\end{Corollary}

The explicit characterization $\Iaux_{\sf Cartan} \equiv 0$ 
of sphericity for {\em tubes} appeared 
{\em e.g.} in Loboda's 
articles~{\cite{Loboda-1996, Loboda-1997, Loboda-2011}},
{\em cf.} also~{\cite{Isaev-2011}} in $\C^{n+1}$ for
any $n \geqslant 1$.

In this much studied tube context, {\em affine} equivalence
to the model parabola $\{ u = x^2 \}$ is
characterized by the vanishing of a {\em different}
invariant. 

\begin{Theorem}
{\rm {\cite{Halphen-1878, Chen-Merker-2020}}}
The following two conditions are equivalent for a $\mathcal{C}^\omega$
curve $\gamma = \big\{ u = F(x) \big\}$ in the plane
$\R_{x,u}^2$ satisfying $F_{xx} \neq 0$.

\smallskip\noindent{\bf (i)}\,
$\gamma$ is affinely equivalent to $\big\{ u' = (x')^2 \big\}$.

\smallskip\noindent{\bf (ii)}\,
The graphing function $F$ satisfies the 5\textsuperscript{th}
order ordinary differential equation:
\[
0
\,\equiv\,
\Iaux_{\Halphensmall}
\,:=\,
3\,F_{xx}\,F_{xxxx}
-
5\,\big(F_{xxx}\big)^2.
\eqno\qed
\]
\end{Theorem}

It is easy to verify by differentiation that:
\[
\Big(
0
\,\equiv\,
\Iaux_{\Halphensmall}
\Big)
\ \ \ \ \ \ \ \ \ \ \ \ \ \ \ \ \ \
\Longrightarrow
\ \ \ \ \ \ \ \ \ \ \ \ \ \ \ \ \ \
\Big(
\Iaux_{\Cartansmall}
\,\equiv\,
0
\Big),
\]
whereas the reverse implication is false.  So as explained in the
introduction, a classification problem arises, solved by Dadok-Yang
under $\mathcal{C}^7$-smoothness assumption.  We `restrict' their
result to the $\mathcal{C}^\omega$ category, and give their
original statement, equivalent to Isaev's reformulation
(in the case $n = 1$) given in 
Section~{\ref{introduction-affine-rigidity}}.

\begin{Theorem}
{\rm {\cite{Dadok-Yang-1985}}}
Any spherical $\mathcal{C}^\omega$ tube hypersurface
$\big\{ u = F(x) \big\} \subset \C^2$ is equivalent
to one of the following:

\smallskip\noindent{\bf (1)}\,
$u = x^2$;

\smallskip\noindent{\bf (2)}\,
$u = e^x$; 

\smallskip\noindent{\bf (3)}\,
$u = \arcsin\, e^x$;

\smallskip\noindent{\bf (4)}\,
$u = \arcsinh\, e^x$.\qed

\end{Theorem}

\Section{\bf Pocchiola's CR invariants $\Waux_0$ and $\Jaux_0$
\\
for $\mathfrak{C}_{2,1}$ hypersurfaces $M^5 \subset \C^3$}
\label{Pocchiola-invariants-C21-M5-C3}
\HEAD{{\ref{Pocchiola-invariants-C21-M5-C3}}.~{\sf
Pocchiola's CR invariants $\Waux_0$ and $\Jaux_0$
for $\mathfrak{C}_{2,1}$ hypersurfaces $M^5 \subset \C^3$}
}{
Jo\"el {\sc Merker}, 
Paris-Saclay University, Orsay, France}

In a series of papers~{\cite{Isaev-Zaitsev-2013, Isaev-2016,
Isaev-2018}} after a research monograph~{\cite{Isaev-2011}}, Isaev
studied zero CR-curvature equations for a special class of CR
submanifolds $M^5 \subset \C^3$, assuming 
$M^5 = S^2 \times i\,
\R^3$ is a {\sl tube}, with $S^2 \subset \R^3$ a surface. Such
a tube assumption `lightens' differential ring computations. Explicit
(relative) differential invariants are easier to reach.

Recall that the class $\mathfrak{C}_{2,1}$ consists of
$2$-nondegenerate constant Levi rank $1$ hypersurfaces $M^5 \subset
\C^3$. According 
to~{\cite{Merker-Pocchiola-Sabzevari-2013-5-CR-II}}, the
assumption of $2$-nondegeneracy excludes the degenerate situation
where $M^5 \cong \C \times N^3$ is (locally) biholomorphic to a
product of $\C$ with a Levi nondegenerate hypersurface $N^3 \subset
\C^2$.

In this paper, coordinates on $\C^3$ will be equally denoted:
\[
(z_1,z_2,w)
\,=\,
\big(
x_1+i\,y_1,\,
x_2+i\,y_2,\,
u+i\,v
\big)
\ \ \ \ \ \ \ \ \ \ \ \ \ \ \ \ \ \
\text{or}
\ \ \ \ \ \ \ \ \ \ \ \ \ \ \ \ \ \
\big(
x+i\,\zeta,\,
y+i\,\eta,\,
u+i\,v
\big).
\]
In order to avoid Analysis of {\sc pde}'s
(but {\em see} Section~{\ref{smoothness-assumption-improvement}}),
all geometric objects will be assumed 
real-analytic ($\mathcal{C}^\omega$) for the moment

The local biholomorphic equivalence for $M \in \mathfrak{C}_{2,1}$,
especially reduction to an $\{e\}$-structure was studied by
Isaev-Zaitsev in~{\cite{Isaev-Zaitsev-2013}} and by Medori-Spiro
in~{\cite{Medori-Spiro-2014, Medori-Spiro-2015}}, in an abstract CR
setting.  Independently, in an embedded setting, Pocchiola, the
author, and Foo~{\cite{Pocchiola-2013, Merker-Pocchiola-2018,
Foo-Merker-2019}} conducted the Cartan method of equivalence,
doing {\em explicit} calculations in terms of a $\mathcal{C}^\omega$
graphing function:
\[
M
\colon\ \ \ \ \
\big\{
(z_1,z_2,w)\in\C^3\colon\,\,
u
=
F(x_1,y_1,x_2,y_2,v)
\big\}.
\]
Also, Nurowski and the author~{\cite{Merker-Nurowski-2020}} 
classified homogeneous models
for systems of PDEs associated to such $M^5 \in \mathfrak{C}_{2,1}$.
This is useful for application to the classification problem.

The recent prepublication~{\cite{Foo-Merker-2019}}
shows that $\sim\,$50 pages of detailed computations
are required until one arrives at Pocchiola's two 
primary differential invariants:
\[
\Waux_0
\,=\,
\Waux_0\big(J_{x_1,y_1,x_2,y_2,v}^5F\big)
\ \ \ \ \ \ \ \ \ \ \ \ \ \ \ \ \ \
\text{and}
\ \ \ \ \ \ \ \ \ \ \ \ \ \ \ \ \ \
\Jaux_0
\,=\,
\Jaux_0\big(J_{x_1,y_1,x_2,y_2,v}^6F\big).
\]
Secondary invariants are covariant derivatives of $\Waux_0$ and
$\Jaux_0$ within the $\{e\}$-structure bundle.

Now, let us be more specific.  In the intrinsic coordinates $(z_1,
z_2, \overline{ z}_1, \overline{ z}_2, v)$ on $M$, two natural
generators of $T^{1,0}M$ are:
\[
\mathcal{L}_1
\,:=\,
\frac{\partial}{
\partial z_1}
-
i\,
\frac{F_{z_1}}{
1+i\,F_v}\,
\frac{\partial}{
\partial v}
\ \ \ \ \ \ \ \ \ \ \ \ \ \
\text{and}
\ \ \ \ \ \ \ \ \ \ \ \ \ \
\mathcal{L}_2
\,:=\,
\frac{\partial}{
\partial z_2}
-
i\,
\frac{F_{z_2}}{
1+i\,F_v}\,
\frac{\partial}{\partial v}.
\]
Then $\overline{ \mathcal{L}}_1$ and $\overline{ \mathcal{L}}_2$ 
generate the conjugate bundle
$T^{0,1}M = \overline{T^{1,0}M}$. Abbreviate:
\[
\Aaux^1
\,:=\,
-\,i\,
\frac{\Faux_{z_1}}{
1+i\,\Faux_v}
\ \ \ \ \ \ \ \ \ \ \ \ \ \ \ \ \ \ \ \
\text{and}
\ \ \ \ \ \ \ \ \ \ \ \ \ \ \ \ \ \ \ \
\Aaux^2
\,:=\,
-\,i\,
\frac{\Faux_{z_2}}{
1+i\,\Faux_v}.
\]

Clearly, the real differential $1$-form:
\[
\varrho_0
\,:=\,
dv
-
\Aaux^1\,dz_1
-
\Aaux^2\,dz_2
-
\overline{\Aaux}^1\,
d\overline{z}_1
-
\overline{\Aaux}^2\,
d\overline{z}_2
\]
has kernel the sum of these two bundles:
\[
\big\{
\varrho_0
=
0
\big\}
\,=\,
T^{1,0}M
\oplus
T^{0,1}M.
\]
At various points $p = \big( z_1, z_2, \overline{ z}_1, \overline{
z}_2, v \big)$ on $M$, and in terms of this $1$-form $\varrho_0$,
the hypothesis that $M$ has everywhere degenerate Levi form 
reads as:
\[
\left\vert\!
\begin{array}{cc}
\varrho_0\big(i\,[\mathcal{L}_1,\overline{\mathcal{L}}_1]\big)
&
\varrho_0\big(i\,[\mathcal{L}_2,\overline{\mathcal{L}}_1]\big)
\\
\varrho_0\big(i\,[\mathcal{L}_1,\overline{\mathcal{L}}_2]\big)
&
\varrho_0\big(i\,[\mathcal{L}_2,\overline{\mathcal{L}}_2]\big)
\end{array}
\!\right\vert
(p)
\,=\,
0
\eqno
{\scriptstyle{(\forall\,p\,\in\,M)}}.
\]

The assumption that the Levi form has 
constant rank equal to $1$
(but not $0$!)
expresses as the fact that the (real) vector field:
\[
\mathcal{T}
\,:=\,
i\,
\big[\mathcal{L}_1,\overline{\mathcal{L}}_1\big]
\,=\,
i\,
\Big(
\mathcal{L}_1\big(\overline{\Aaux}^1\big)
-
\overline{\mathcal{L}}_1\big(\Aaux^1\big)
\Big)
\frac{\partial}{\partial v}
\,=:\,
\laux\,
\frac{\partial}{\partial v},
\]
is nowhere vanishing, {\em i.e.}:
\[
0
\,\neq\,
\laux
\,:=\,
i\,
\Big(\overline{\Aaux}_{z_1}^1
+
\Aaux^1\,\overline{\Aaux}_v^1-
\Aaux_{\overline{z}_1}^1
-\overline{\Aaux}^1\,\Aaux_v^1\Big).
\]

The Levi kernel bundle $K^{1,0}M \subset T^{1,0}M$
is then generated by the $(1,0)$-vector field:
\[
\mathcal{K}
\,:=\,
\kaux\,\mathcal{L}_1
+
\mathcal{L}_2,
\]
with the fundamental {\sl slant function}:
\[
\kaux
\,:=\,
-\,
\frac{
\mathcal{L}_2\big(\overline{\Aaux}^1\big)
-
\overline{\mathcal{L}}_1\big(\Aaux^2\big)}{
\mathcal{L}_1\big(\overline{\Aaux}^1\big)
-
\overline{\mathcal{L}}_1\big(\Aaux^1\big)}.
\]
The assumption of $2$-nondegeneracy is then equivalent 
{\cite{Merker-Pocchiola-Sabzevari-2013-5-CR-II,
Pocchiola-2013, Merker-Pocchiola-2018}} 
to the nonvanishing:
\[
0
\,\neq\,
\overline{\mathcal{L}}_1(\kaux).
\]
Also, the conjugate field $\overline{\mathcal{K}}$ generates
the conjugate Levi kernel bundle
$K^{0,1}M \subset T^{0,1}M$. 

Similary as for hypersurfaces $M^3 \subset \C^2$, there also is a
second fundamental function:
\[
\Paux
\,:=\,
\frac{\laux_{z_1}+\Aaux^1\,\laux_v-\laux\,\Aaux_v^1}{\laux}.
\]

Next, introduce the five $1$-forms:
\[
\aligned
\rho_0
&
=
\frac{dv-\Aaux^1dz_1-\Aaux^2dz_2
-\overline{\Aaux}^1d\overline{z}_1
-\overline{\Aaux}^2d\overline{z}_2}{\laux},
\\
\kappa_0
&
=
dz_1-\kaux\,dz_2,
\\
\zeta_0
&
=
dz_2,
\\
\overline{\kappa}_0
&
=
d\overline{z}_1
-
\overline{\kaux}\,d\overline{z}_2,
\\
\overline{\zeta}_0
&
=
d\overline{z}_2.
\endaligned
\]
After intensive computations, redone in~{\cite{Foo-Merker-2019}},
Pocchiola obtained modifications $\big\{ \rho, \kappa, \zeta,
\overline{ \kappa}, \overline{ \zeta} \big\}$ of these $1$-forms
$\big\{ \rho_0, \kappa_0, \zeta_0, \overline{\kappa}_0,
\overline{\zeta}_0 \big\}$, together with certain $1$-forms $\pi^1$,
$\pi^2$, $\overline{ \pi}^1$, $\overline{\pi}^2$ which satisfy
structure equations of the specific concise shape:
\leqnomode\usetagform{default}
\begin{align}
\label{concise-final-drho-dkappa-dzeta}
d\rho
&
\,=\,
\big(
\pi^1
+
\overline{\pi}^1
\big)
\wedge\rho
+
i\,\kappa\wedge\overline{\kappa},
\notag
\\
d\kappa
&
\,=\,
\pi^2\wedge\rho
+
\pi^1\wedge\kappa
+
\zeta\wedge\overline{\kappa},
\notag
\\
d\zeta
&
\,=\,
\big(\pi^1-\overline{\pi}^1\big)
\wedge\zeta
+
i\,\pi^2\wedge\kappa
\,+
\\
&
\ \ \ \ \
+
\Raux\,
\rho\wedge\zeta
+
i\,
\frac{1}{\overline{\sf c}^3}\,
\overline{\Jaux}_0\,\rho\wedge\overline{\kappa}
+
\frac{1}{{\sf c}}\,
\Waux_0\,
\kappa\wedge\zeta,
\notag
\end{align}
in which $\Raux$ is a secondary invariant:
\[
\Raux
\,:=\,
\Re\,
\left[
i\,
\frac{{\sf e}}{{\sf c}{\sf c}}\,
\Waux_0
+
\frac{1}{{\sf c}\overline{\sf c}}
\bigg(
-\,\frac{i}{2}\,
\overline{\mathcal{L}}_1\big(\Waux_0\big)
+
\frac{i}{2}\,
\bigg(
-\,\frac{1}{3}\,
\frac{\overline{\mathcal{L}}_1\big(
\overline{\mathcal{L}}_1(\kaux)\big)}{
\overline{\mathcal{L}}_1(\kaux)}
+
\frac{1}{3}\,
\overline{\Paux}
\bigg)\,
\Waux_0
\bigg)
\right],
\]
expressed in terms of
Pocchiola's two primary invariants:
\[
\aligned
\Waux_0
&
\,:=\,
-\,\frac{1}{3}\,
\frac{\mathcal{K}\big(\overline{\mathcal{L}}_1\big(
\overline{\mathcal{L}}_1(\kaux)\big)\big)}{
\overline{\mathcal{L}}_1(\kaux)^2}
+
\frac{1}{3}\,
\frac{\mathcal{K}\big(\overline{\mathcal{L}}_1(\kaux)\big)\,\,
\overline{\mathcal{L}}_1\big(\overline{\mathcal{L}}_1(\kaux)\big)}{
\overline{\mathcal{L}}_1(\kaux)^3}
\,+
\\
&
\ \ \ \ \
+
\frac{2}{3}\,
\frac{\mathcal{L}_1\big(\mathcal{L}_1(\overline{\kaux})\big)}{
\mathcal{L}_1(\overline{\kaux})}
+
\frac{2}{3}\,
\frac{\mathcal{L}_1\big(\overline{\mathcal{L}}_1(\kaux)\big)}{
\overline{\mathcal{L}}_1(\kaux)}
+
\frac{i}{3}\,
\frac{\mathcal{T}(\kaux)}{\overline{\mathcal{L}}_1(\kaux)},
\\
\overline{\Jaux}_0
&
\,:=\,
\frac{1}{6}\,
\frac{\overline{\mathcal{L}}_1\big(
\overline{\mathcal{L}}_1\big(
\overline{\mathcal{L}}_1\big(
\overline{\mathcal{L}}_1(\kaux)\big)\big)\big)}{
\overline{\mathcal{L}}_1(\kaux)}
-
\frac{5}{6}\,
\frac{\overline{\mathcal{L}}_1\big(
\overline{\mathcal{L}}_1\big(
\overline{\mathcal{L}}_1(\kaux)\big)\big)\,\,
\overline{\mathcal{L}}_1\big(
\overline{\mathcal{L}}_1(\kaux)\big)
}{
\overline{\mathcal{L}}_1(\kaux)^2}
-
\frac{1}{6}\,
\frac{\overline{\mathcal{L}}_1\big(
\overline{\mathcal{L}}_1\big(
\overline{\mathcal{L}}_1(\kaux)\big)\big)
}{
\overline{\mathcal{L}}_1(\kaux)}\,
\overline{\Paux}
\,+
\\
&
\ \ \ \ \
+
\frac{20}{27}\,
\frac{\overline{\mathcal{L}}_1\big(\overline{\mathcal{L}}_1
(\kaux)\big)^3}{
\overline{\mathcal{L}}_1(\kaux)^3}
+
\frac{5}{18}\,
\frac{\overline{\mathcal{L}}_1\big(
\overline{\mathcal{L}}_1(\kaux)\big)^2}{
\overline{\mathcal{L}}_1(\kaux)^2}\,
\overline{\Paux}
+
\frac{1}{6}\,
\frac{\overline{\mathcal{L}}_1\big(
\overline{\mathcal{L}}_1(\kaux)\big)\,\,
\overline{\mathcal{L}}_1\big(\overline{\Paux}\big)}{
\overline{\mathcal{L}}_1(\kaux)}
-
\frac{1}{9}\,
\frac{\overline{\mathcal{L}}_1\big(
\overline{\mathcal{L}}_1(\kaux)\big)}{
\overline{\mathcal{L}}_1(\kaux)}\,\,
\overline{\Paux}\,\overline{\Paux}
\,-
\\
&
\ \ \ \ \
-
\frac{1}{6}\,
\overline{\mathcal{L}}_1\big(
\overline{\mathcal{L}}_1\big(
\overline{\Paux}\big)\big)
+
\frac{1}{3}\,
\overline{\mathcal{L}}_1\big(\overline{\Paux}\big)\,
\overline{\Paux}
-
\frac{2}{27}\,
\overline{\Paux}\,
\overline{\Paux}\,
\overline{\Paux}.
\endaligned
\] 
In depth and quite strikingly, the numerators of $\Waux$ and $\Jaux$
both contain $> {\bf 10^4}$ differential jet monomials.
Fortunately, when $M$ is assumed to be tube, 
we will soon see how simpler $\Waux_0$ and $\Jaux_0$ become.

Without any special assumption on $F$, a byproduct of Cartan's method
characterizes hypersurfaces $M^5 \subset \C^3$ having zero Pocchiola
curvature, as being biholomorphically equivalent to a known {\em
model}.

\begin{Theorem}
{\rm \cite{Pocchiola-2013, Merker-Pocchiola-2018, Foo-Merker-2019}}
For a $\mathcal{C}^\omega$ hypersurface
$M^5 \subset \C^3$ belonging to the class $\mathfrak{C}_{2,1}$,
the following two conditions are equivalent:

\smallskip\noindent{\bf (i)}\,
$0 \equiv \Waux_0 \equiv \Jaux_0$;

\smallskip\noindent{\bf (ii)}\,
$M^5 \subset \C^3$ is locally biholomorphic to the CR tube:
\[
T
\,:=\,
\Big\{
(z,\zeta,w)
\in
\C^3
\colon\,\,\,
\Re\,w
\,=\,
\frac{(\Resmall\,z)^2}{1-\Resmall\,\zeta}
\Big\}.
\eqno\qed
\]
\end{Theorem}

To state the {\sl relative invariancy} property satisfied by $\Waux_0$
and $\Jaux_0$, let us introduce the

\begin{Notation}
The symbol `$\nonzero$' shall denote various local
$\mathcal{C}^\omega$ or $\mathcal{C}^\infty$ functions which are {\em
nowhere vanishing}\,\,---\,\,possibly after restriction to some 
smaller open subsets.
\end{Notation}

Indeed, general Cartan method guarantees
that $\Waux_0$ and $\Jaux_0$ are {\sl relative invariants}
in the following sense.
Suppose $h \colon \C^3 \longrightarrow {\C'}^3$ is a local
biholomorphism which sends CR-diffeomorphically 
$M$ onto its image $M' := h(M)$, graphed similarly as:
\[
M'
\colon\ \ \ \ \
\big\{
(z_1',z_2',w')\in\C^3\colon\,\,
u'
=
F'\big(x_1',y_1',x_2',y_2',v'\big)
\big\}.
\]
Pocchiola's invariants for $M'$ are computed by means
of {\em exactly the same} universal formulas
in terms of $F'$.

\begin{Theorem}
Under a biholomorphic equivalence:
\[
\Waux_0\big(F'\big)
\,=\,
\nonzero
\cdot
\Waux_0\big(F\big)
\ \ \ \ \ \ \ \ \ \ \ \ \ \ \ \ \ \
\text{and}
\ \ \ \ \ \ \ \ \ \ \ \ \ \ \ \ \ \
\Jaux_0\big(F'\big)
\,=\,
\nonzero
\cdot
\Jaux_0\big(F\big).
\eqno\qed
\]
\end{Theorem}

As an obvious corollary:
\[
\big\{
0
=
\Waux_0(F')
=
\Jaux_0(F')
\big\}
\,\,=\,\,
\big\{
0
=
\Waux_0(F)
=
\Jaux_0(F)
\big\}.
\]

Next, let us come back to the affine transformation group $\Aff(\R^3)$
presented in Section~{\ref{introduction-affine-rigidity}}.  It is clear
that all {\em real} affine transformations $\xaux \longmapsto A \xaux
+ b$ of $\R^3$ extend as biholomorphic
transformations $\zaux \longmapsto A \zaux + b$ of $\C^3$, with $\zaux
= (z, \zeta, w)$ and $\xaux = \Re\, \zaux$. Although the group
inclusion:
\[
\Aff(\R^3)
\,\,\subset\,\,
\Bihol(\C^3),
\]
shows a high-dimensional discrepancy:
\[
12
\,\,<\,\,
\infty,
\] 
we may deduce by pure `logic' that the `{\sl affinizations}' of
Pocchiola's invariants:
\[
\Waux_{\sf aff}
\,:=\,
\Waux_0
\big\vert_{
M\,\text{is tube}},
\ \ \ \ \ \ \ \ \ \ \ \ \ \ \ \ \ \ \ \ \ \ \ \ \ \
\Jaux_{\sf aff}
\,:=\,
\Jaux_0
\big\vert_{
M\,\text{is tube}},
\]
namely:
\[
\Waux_{\sf aff}
\,=\,
\Waux_0
\big(
F(x_1,x_2)
\big),
\ \ \ \ \ \ \ \ \ \ \ \ \ \ \ \ \ \ \ \ \ \ \ \ \ \
\Jaux_{\sf aff}
\,=\,
\Jaux_0
\big(
F(x_1,x_2)
\big),
\] 
are {\em also} relative invariants under {\sl affine} transformations
of $\R^3$.  More precisely, if
$g \colon \R^3 \longrightarrow \R^3$
denotes any affine (invertible) map which sends a surface $S = \big\{
u = F(x_1, x_2) \big\}$ onto its image $S' := g(S)$, graphed similarly
as $\big\{ u' = F'(x_1', x_2') \big\}$, we deduce

\begin{Theorem}
Under a real affine equivalence of $\R^3$:
\[
\Waux_{\sf aff}\big(F'\big)
\,=\,
\nonzero
\cdot
\Waux_{\sf aff}\big(F\big)
\ \ \ \ \ \
\text{and}
\ \ \ \ \ \
\Jaux_{\sf aff}\big(F'\big)
\,=\,
\nonzero
\cdot
\Jaux_{\sf aff}\big(F\big).
\eqno\qed
\]
\end{Theorem}

The next Section~{\ref{W-aff-J-aff}} is devoted to show how to clean
up appropriate explicit expressions for $\Waux_{\sf aff}$ and
$\Jaux_{\sf aff}$. Also,
Section~{\ref{affine-invariants-graph-transforms}} endeavors to
recover from scratch the (relative) invariancy of $\Waux_{\sf aff}$
and of $\Jaux_{\sf aff}$ under affine transformations of the real
space $\R^3 \ni (x,y,u)$.

Next, because $\dim\, \Bihol_3(\C) \gg \dim\, \Aff_3(\R)$, 
and in view of the $\C^2$ case presented in
Section~{\ref{prologue-C2-case}}, it is
natural to expect that there exist hypersurfaces $M^5 \in
\mathfrak{C}_{2,1}$ such that:
\[
M
\overset{\Biholsmall}{\,\,\cong\,\,}
T,
\ \ \ \ \ \ \ \ \ \ \ \ \ \ \ \ \ \
\text{while}
\ \ \ \ \ \ \ \ \ \ \ \ \ \ \ \ \ \
M
\overset{\Affsmall}{\,\,\not\cong\,\,}
T.
\]

One could try to 
find affine (relative)
differential invariants $\Iaux_1, \Iaux_2, \dots$,
whose vanishing characterizes affine equivalence of a surface $S =
\big\{ u = F(x,y) \big\}$ to the model:
\[
0
\,\equiv\,
\Iaux_1
\,\equiv\,
\Iaux_2
\,\equiv\,
\cdots
\ \ \ \ \ \ \ \ \ \ \ \ \ \ \ \ \ \
\Longleftrightarrow
\ \ \ \ \ \ \ \ \ \ \ \ \ \ \ \ \ \
S
\overset{\Affsmall}{\,\,\cong\,\,}
\Big\{
u
\,=\,
\frac{x^2}{1-y}
\Big\}.
\]

Of course, $\Waux_{\sf aff}$ and $\Jaux_{\sf aff}$ 
are among $\Iaux_1, \Iaux_2, \dots$. 
So the question is: are there further affine
invariants? It might very well be so!

\Section{\bf Affine Pocchiola Invariants $\Waux_{\sf aff}$ and
$\Jaux_{\sf aff}$ 
\\
for Tube hypersurfaces $M^5 = S^2 \times 
\big( i\, \R^3 \big) \subset \C^3$}
\label{W-aff-J-aff}
\HEAD{{\ref{W-aff-J-aff}}.~{\sf Affine Pocchiola Invariants 
$\Waux_{\sf aff}$ and
$\Jaux_{\sf aff}$ for Tube hypersurfaces $M^5 = S^2 \times 
\big( i\, \R^3 \big) \subset \C^3$}
}{
Jo\"el {\sc Merker}, 
Paris-Saclay University, Orsay, France}

Suppose therefore that $M^5 = S^2 \times i\, \R^3$
is tube:
\[
\big\{
u
\,=\,
F(x_1,x_2)
\big\}.
\]
Then:
\[
\aligned
\mathcal{L}_1
&
\,=\,
\frac{\partial}{\partial z_1}
-
\frac{i}{2}\,F_{x_1}\,\frac{\partial}{\partial v},
\ \ \ \ \ \ \ \ \ \ \ \ \ \ \ \ \ \ \ \ \ \ \ \ \ \
\mathcal{L}_2
\,=\,
\frac{\partial}{\partial z_2}
-
\frac{i}{2}\,F_{x_2}\,\frac{\partial}{\partial v},
\\
\overline{\mathcal{L}}_1
&
\,=\,
\frac{\partial}{\partial\overline{z}_1}
+
\frac{i}{2}\,F_{x_1}\,\frac{\partial}{\partial v},
\ \ \ \ \ \ \ \ \ \ \ \ \ \ \ \ \ \ \ \ \ \ \ \ \ \
\overline{\mathcal{L}}_2
\,=\,
\frac{\partial}{\partial\overline{z}_2}
+
\frac{i}{2}\,F_{x_2}\,\frac{\partial}{\partial v},
\endaligned
\]
whence:
\[
\mathcal{K}
\,=\,
\kaux\,\mathcal{L}_1
+
\mathcal{L}_2
\,=\,
-\,\frac{F_{x_1x_2}}{F_{x_1x_1}}\,
\mathcal{L}_1
+
\mathcal{L}_2.
\]
So the action of the derivations $\mathcal{L}_1$, $\mathcal{K}$,
$\overline{\mathcal{L}}_1$, $\overline{\mathcal{K}}$
on functions depending only on 
$(x_1, x_2)$ identifies with the actions of the
{\em purely real} vector fields:
\[
\aligned
L_1
&
\,:=\,
\frac{1}{2}\,
\frac{\partial}{\partial x_1},
\\
K
&
\,:=\,
-\,\frac{1}{2}\,
\frac{F_{x_1x_2}}{F_{x_1x_1}}\,
\frac{\partial}{\partial x_1}
+
\frac{1}{2}\,
\frac{\partial}{\partial x_2}.
\endaligned
\]
It follows that all four quantities:
\[
\overline{\mathcal{L}}_1(\kaux)
\,=\,
\mathcal{L}_1(\overline{\kaux})
\,=\,
\mathcal{L}_1(\kaux)
\,=\,
\overline{\mathcal{L}}_1(\overline{\kaux})
\,=\,
-\,\frac{1}{2}\,
\frac{F_{xx}\,F_{xxy}-F_{xy}\,F_{xxx}}{(F_{xx})^2}
\]
are real, where we already have switched notation:
\[
(x_1,x_2)
\,\equiv\,
(x,y).
\]
Then the second fundamental function is also real:
\[
\Paux
\,=\,
\frac{1}{2}\,
\frac{F_{xxx}}{F_{xx}}
\,=\,
\overline{\Paux}.
\]

Observe from reality the vanishing:
\[
\mathcal{T}(\kaux)
\,=\,
i\,\big[\mathcal{L}_1,\overline{\mathcal{L}}_1\big](\kaux)
\,=\,
i\,\mathcal{L}_1\big(\overline{\mathcal{L}}_1(\kaux)\big)
-
i\,\overline{\mathcal{L}}_1\big(\mathcal{L}_1(\kaux)\big)
\,\,=\,\,
0.
\]
By reading and translating $\Waux_0$ and $\Jaux_0$ above, we obtain:
\[
\aligned
\Waux_{\sf aff}
&
\,=\,
\frac{2}{3}\,
\frac{L_1\big(L_1(\kaux))}{L_1(\kaux)}
+
\frac{2}{3}\,
\frac{L_1\big(L_1(\kaux))}{L_1(\kaux)}
\,+
\\
&
\ \ \ \ \
+
\frac{1}{3}\,
\frac{L_1\big(L_1(\kaux)\big)\,K\big(L_1(\kaux)\big)}{
L_1(\kaux)^3}
-
\frac{1}{3}\,
\frac{K\big(L_1\big(L_1(\kaux)\big)\big)}{L_1(\kaux)^2}
+
0,
\endaligned
\]
together with:
\[
\aligned
\Jaux_{\sf aff}
&
\,=\,
\frac{1}{6}\,
\frac{L_1\big(L_1\big(L_1\big(L_1(\kaux)\big)\big)\big)}{
L_1(\kaux)}
-
\frac{5}{6}\,
\frac{L_1\big(L_1\big(L_1(\kaux)\big)\big)\,\,
L_1\big(L_1(\kaux)\big)}{L_1(\kaux)^2}
-
\frac{1}{6}\,
\frac{L_1\big(L_1\big(L_1(\kaux)\big)\big)}{L_1(\kaux)}\,\Paux
\,+
\\
&
\ \ \ \ \
+
\frac{20}{27}\,
\frac{L_1\big(L_1(\kaux)\big)}{L_1(\kaux)^3}
+
\frac{5}{18}\,
\frac{L_1\big(L_1(\kaux)\big)^2}{L_1(\kaux)}\,\Paux
+
\frac{1}{6}\,
\frac{L_1\big(L_1(\kaux)\big)\,L_1(\Paux)}{L_1(\kaux)}
-
\frac{1}{9}\,
\frac{L_1\big(L_1(\kaux)\big)}{L_1(\kaux)}{\Paux\,\Paux}
\,-
\\
&
\ \ \ \ \
-\,
\frac{1}{6}\,
L_1\big(L_1(\kaux)\big)
+
\frac{1}{3}\,
L_1(\Paux)\,\Paux
-
\frac{2}{27}\,
\Paux\,\Paux\,\Paux.
\endaligned
\]

The expansion of $\Jaux_{\sf aff}$ can be done plainly:
\[
\footnotesize
\aligned
\Jaux_{\sf aff}
&
\,:=\,
-\,
\frac{1}{54}\,
\frac{1}{\big(
F_{xxy}\,F_{xx}-F_{xy}\,F_{xxx}
\big)^3}\,
\bigg\{
\\
&
\ \ \ \ \ \ \
-\,9\,F_{xxxxxy}\,F_{xx}^3\,F_{xxy}^2
+
45\,F_{xxxxx}\,F_{xx}^2\,F_{xxy}^3
-
45\,F_{xy}^2\,F_{xxx}^3\,F_{xxxxy}
+
\\
&
\ \ \ \ \ \ \
+
9\,F_{xy}^3\,F_{xxxxxx}\,F_{xxx}^2
-
40\,F_{xxxy}^3\,F_{xx}^3
+
40\,F_{xy}^3\,F_{xxxx}^3
-
90\,F_{xxx}\,F_{xy}^2\,F_{xxxx}^2\,F_{xxy}
+
\\
&
\ \ \ \ \ \ \
+
45\,F_{xy}^2\,F_{xxx}^2\,F_{xxxxx}\,F_{xxy}
-
45\,F_{xy}^3\,F_{xxxx}\,F_{xxxxx}\,F_{xxx}
+
90\,F_{xxxy}\,F_{xy}^2\,F_{xxxx}\,F_{xxx}^2
+
\\
&
\ \ \ \ \ \ \
+
90\,F_{xxx}\,F_{xxxy}^2\,F_{xx}^2\,F_{xxy}
-
90\,F_{xxxy}^2\,F_{xx}\,F_{xy}\,F_{xxx}^2
+
120\,F_{xxxy}^2\,F_{xx}^2\,F_{xy}\,F_{xxxx}
-
\\
&
\ \ \ \ \ \ \
-\,
120\,F_{xxxy}\,F_{xx}\,F_{xy}^2\,F_{xxxx}^2
-
90\,F_{xxxy}\,F_{xx}^2\,F_{xxy}^2\,F_{xxxx}
-
45\,F_{xxy}^2\,F_{xxx}\,F_{xx}^2\,F_{xxxxy}
+
\\
&
\ \ \ \ \ \ \
+
90\,F_{xy}\,F_{xxxx}^2\,F_{xx}\,F_{xxy}^2
+
45\,F_{xxxy}\,F_{xx}^3\,F_{xxxxy}\,F_{xxy}
-
9\,F_{xxxxxy}\,F_{xx}\,F_{xy}^2\,F_{xxx}^2
+
\\
&
\ \ \ \ \ \ \
+
9\,F_{xy}\,F_{xxxxxx}\,F_{xx}^2\,F_{xxy}^2
-
45\,F_{xxxy}\,F_{xx}^2\,F_{xy}\,F_{xxxxx}\,F_{xxy}
+
45\,F_{xxxy}\,F_{xx}\,F_{xy}^2\,F_{xxxxx}\,F_{xxx}
+
\\
&
\ \ \ \ \ \ \
+
90\,F_{xxy}\,F_{xxx}^2\,F_{xx}\,F_{xxxxy}\,F_{xy}
-
45\,F_{xy}\,F_{xxxx}\,F_{xx}^2\,F_{xxxxy}\,F_{xxy}
+
45\,F_{xy}^2\,F_{xxxx}\,F_{xx}\,F_{xxxxy}\,F_{xxx}
+
\\
&
\ \ \ \ \ \ \
+
45\,F_{xy}^2\,F_{xxxx}\,F_{xx}\,F_{xxxxx}\,F_{xxy}
-
45\,F_{xxxy}\,F_{xx}^2\,F_{xxxxy}\,F_{xy}\,F_{xxx}
-
90\,F_{xxy}^2\,F_{xxx}\,F_{xx}\,F_{xy}\,F_{xxxxx}
+
\\
&
\ \ \ \ \ \ \
+
18\,F_{xxxxxy}\,F_{xx}^2\,F_{xy}\,F_{xxx}\,F_{xxy}
-
18\,F_{xy}^2\,F_{xxxxxx}\,F_{xx}\,F_{xxx}\,F_{xxy}
\bigg\}.
\endaligned
\]
However, in the expansion of $\Waux_{\sf aff}$, one must take
account of relations coming from the assumption
that the real Hessian of $F$ vanishes identically:
\[
F_{yy}
\,=\,
\frac{\big(F_{xy}\big)^2}{F_{xx}}.
\]
Differentiations with respect to $x$ and to $y$ followed by
replacements give:
\[
\aligned
F_{xyy}
&
\,=\,
2\,
\frac{F_{xy}\,F_{xxy}}{F_{xx}}
-
\frac{\big(F_{xy}\big)^2\,F_{xxx}}{\big(F_{xx}\big)^2},
\\
F_{yyy}
&
\,=\,
3\,
\frac{\big(F_{xy}\big)^2\,F_{xxy}}{\big(F_{xx}\big)^2}
-
2\,
\frac{\big(F_{xy}\big)^3\,F_{xxx}}{\big(F_{xx}\big)^3}.
\endaligned
\]
Next:
\[
\footnotesize
\aligned
F_{xxyy}
&
\,=\,
2\,\frac{\big(F_{xxy}\big)^2}{F_{xx}}
-
4\,\frac{F_{xy}\,F_{xxy}\,F_{xxx}}{\big(F_{xx}\big)^2}
+
2\,\frac{F_{xy}\,F_{xxxy}}{F_{xx}}
+
2\,\frac{\big(F_{xy}\big)^2\,\big(F_{xxx}\big)^2}{\big(F_{xx}\big)^3}
-
\frac{\big(F_{xy}\big)^2\,F_{xxxx}}{\big(F_{xx}\big)^2},
\\
F_{xyyy}
&
\,=\,
6\,\frac{F_{xy}\,\big(F_{xxy}\big)^2}{\big(F_{xx}\big)^2}
-
12\,\frac{\big(F_{xy}\big)^2\,F_{xxx}\,F_{xxy}}{\big(F_{xx}\big)^3}
+
3\,\frac{\big(F_{xy}\big)^2\,F_{xxxy}}{\big(F_{xx}\big)^2}
+
6\,\frac{\big(F_{xy}\big)^3\,\big(F_{xxx}\big)^2}{\big(F_{xx}\big)^4}
-
2\,\frac{\big(F_{xy}\big)^3\,F_{xxxx}}{\big(F_{xx}\big)^3},
\\
F_{yyyy}
&
\,=\,
12\,\frac{\big(F_{xy}\big)^2\,\big(F_{xxy}\big)^2}{\big(F_{xx}\big)^3}
-
24\,\frac{\big(F_{xy}\big)^3\,F_{xxx}\,F_{xxy}}{\big(F_{xx}\big)^4}
+
12\,\frac{\big(F_{xy}\big)^4\,\big(F_{xxx}\big)^2}{\big(F_{xx}\big)^5}
+
4\,\frac{\big(F_{xy}\big)^3\,F_{xxxy}}{\big(F_{xx}\big)^3}
-
3\,\frac{\big(F_{xy}\big)^4\,F_{xxxx}}{\big(F_{xx}\big)^4}.
\endaligned
\]
Similar formulas exist for $F_{xxxyy}$, $F_{xxyyy}$, $F_{xyyyy}$,
$F_{yyyyy}$.

With a different approach, 
Isaev found in~{\cite{Isaev-2016, Isaev-2018}}
that after these replacements, $\Waux_{\sf aff}$ 
which seems to be a 5\textsuperscript{th}-order
invariant, is in fact a 4\textsuperscript{th}-order one.

\begin{Proposition}
After plain replacements:
\[
\Waux_{\sf aff}
\,=\,
\frac{\big(F_{xx}\big)^2\,F_{xxxy}
-
F_{xx}\,F_{xy}\,F_{xxxx}
+
2\,F_{xy}\,\big(F_{xxx}\big)^2
-
2\,F_{xx}\,F_{xxx}\,F_{xxy}}{
F_{xx}\,\big(F_{xx}\,F_{xxy}-F_{xy}\,F_{xxx}\big)^2}.
\eqno\qed
\]
\end{Proposition}

Then under the hypothesis $0 \equiv \Waux_{\sf aff}$,
many terms in $\Jaux_{\sf aff}$ above cancel,
and if we denote:
\[
\aligned
\Jaux_{\sf aff}^\sim
\,:=\,
&\,
\Jaux_{\sf aff}
\ \ \
\mod\,\,\Waux_{\sf aff}
\\
\,=\,
&\,
-\,
\frac{1}{6}\,
L_1\big(L_1(\kaux)\big)
+
\frac{1}{3}\,
L_1(\Paux)\,\Paux
-
\frac{2}{27}\,
\Paux\,\Paux\,\Paux
\ \ \ \ \
\Big(
\text{\rm using}\,\,
\Waux_{\sf aff}
\,\equiv\,
0
\Big),
\endaligned
\]
this object simplifies as:
\[
\Jaux_{\sf aff}^\sim
\,=\,
-\,\frac{1}{432}\,
\frac{
9\,\big(F_{xx}\big)^2\,
F_{xxxxx}
-
45\,F_{xx}\,F_{xxx}\,F_{xxxx}
+
40\,\big(F_{xxx}\big)^3}{
\big(F_{xx}\big)^3}.
\]
We recognize the Monge invariant with respect to the
first variable $x$, whose vanishing characterizes
the fact that a planar graphed curve $\big\{ u = F(x) \big\}$
in $\R_{x,u}^2$ is contained in a (nondegenerate) conic
({\cite{Halphen-1878, Chen-Merker-2020}} and {\em see} also 
Section~{\ref{affine-invariants-graph-transforms}}).

However, one may convince oneself that $\Jaux_{\sf aff}^\sim =
\Jaux_{\sf aff}\,\, \mod\, \Waux_{\sf aff}$ is {\em not} an affine
relative invariant. Anyway, the two zero-sets coincide:
\[
\big\{ 
0 
\equiv 
\Waux_{\sf aff} 
\equiv 
\Jaux_{\sf aff} 
\big\}
\,\,=\,\,
\big\{ 
0 
\equiv 
\Waux_{\sf aff} 
\equiv 
\Jaux_{\sf aff}^\sim
\big\},
\]
and in conclusion, we may formulate a

\begin{Proposition}
{\rm {\cite{Isaev-2018}}}
CR-flatness of hypersurfaces $M \in \mathfrak{C}_{2,1}$ 
that are tube $\big\{ u = F (x,y) \big\}$ is characterized by:
\reqnomode\usetagform{EngelLie}
\begin{align}
0
&
\,\equiv\,
2\,F_{xy}\,
\big(F_{xxx}\big)^2
-
2\,F_{xx}\,
F_{xxx}\,
F_{xxy}
+
\big(
F_{xx}
\big)^2\,
F_{xxxy}
-
F_{xx}\,F_{xy}\,
F_{xxxx},
\notag
\\
0
&
\,\equiv\,
9\,\big(F_{xx}\big)^2
F_{xxxxx}
-
45\,F_{xx}\,F_{xxx}\,F_{xxxx}
+
40\,\big(F_{xxx}\big)^3.
\tag{\qed}
\end{align}
\end{Proposition}

Once these equations have been obtained
and cleaned up, we can now present
our short proof of Theorem~{\ref{Thm-main}},
in the $\mathcal{C}^\omega$ category.

\Section{\bf Affine Rigidity via Differential Algebra Elimination}
\label{affine-rigidity-elimination}
\HEAD{{\ref{affine-rigidity-elimination}}.~{\sf Affine Rigidity 
via Differential Algebra Elimination}
}{
Jo\"el {\sc Merker}, 
Paris-Saclay University, Orsay, France}

In $\C^3$ with coordinates $\big(z, \zeta, w\big)$, 
with $x = \Re\, z$, $y = \Re\, \zeta$, $u = \Re\, w$,
consider therefore 
a local $\mathcal{C}^\omega$
tube hypersurface graphed as:
\[
M
\colon
\ \ \ \ \
u
\,=\,
F(x,y),
\]
which is of constant Levi rank $1$ and 
$2$-nondegenerate:
\[
F_{xx}
\,\neq\,
0
\,\,\equiv\,\,
F_{xx}\,F_{yy}
-
\big(
F_{xy}
\big)^2
\ \ \ \ \ \ \ \ \ \ \ \ \ \ \ \ \ \
\text{and}
\ \ \ \ \ \ \ \ \ \ \ \ \ \ \ \ \ \
F_{xx}\,
F_{xxy}
-
F_{xy}\,
F_{xxx}
\,\neq\,
0.
\]
As explained 
in~{\cite[Section~2]{Chen-Foo-Merker-Ta-2020}}
the maximally homogeneous
model has an appropriate representation as
the tube:
\[
T_{\sf LC}
\colon
\ \ \ \ \
u
\,=\,
\frac{x^2}{1-y}.
\]

\begin{Theorem}
\label{Theorem-LC-J-W-zero}
A local $\mathcal{C}^\omega$ real surface in $\R^3$:
\[
u
\,=\,
F(x,y)
\]
with $F_{xx} \neq 0$ which has identically zero Hessian: 
\[
0
\,\equiv\,
\Haux_{\sf aff}
\,=\,
\overset{\text{\ding{192}}}{\,\,\equiv\,\,}
F_{xx}\,F_{yy}- 
\big(F_{xy}\big)^2,
\]
is locally affinely equivalent to the model $u =
\frac{x^2}{1-y}$ if and only if:
\[
\aligned
0
&
\overset{\text{\ding{193}}}{\,\,\equiv\,\,}
2\,F_{xy}\,
\big(F_{xxx}\big)^2
-
2\,F_{xx}\,
F_{xxx}\,
F_{xxy}
+
\big(
F_{xx}
\big)^2\,
F_{xxxy}
-
F_{xx}\,F_{xy}\,
F_{xxxx},
\\
0
&
\overset{\text{\ding{194}}}{\,\,\equiv\,\,}
9\,\big(F_{xx}\big)^2
F_{xxxxx}
-
45\,F_{xx}\,F_{xxx}\,F_{xxxx}
+
40\,\big(F_{xxx}\big)^3.
\endaligned
\]
\end{Theorem}

Our (elementary) arguments will consist in normalizing progressively
$F(x,y)$ by means of successive appropriate changes of
affine coordinates, and to `kill' almost all Taylor 
coefficients, thanks to the $3$ equations:
\[
0
\overset{\text{\ding{192}}}{\,\,\equiv\,\,}
\Haux_{\sf aff}
\overset{\text{\ding{193}}}{\,\,\equiv\,\,}
\Waux_{\sf aff}
\overset{\text{\ding{194}}}{\,\,\equiv\,\,}
\Jaux_{\sf aff}^{\sim}.
\]
No integration of any differential equation will be required.

Thus, we recover a result proved by Isaev in~{\cite{Isaev-2016}}.

\begin{Corollary} 
$M\,\,
\text{\em is biholomorphically equivalent to}\,\,
T_{\sf LC}$

\ \ \ \ \ \ \ 
$\,\,\,\Longleftrightarrow\,\,\,
M\,\,
\text{\em is real affinely equivalent to}\,\,
T_{\sf LC}$.\qed
\end{Corollary}

\proof[Proof of Theorem~{\ref{Theorem-LC-J-W-zero}}]
Setting $F(x,y) := \frac{x^2}{1-y}$ in
$\Haux_{\sf aff}$, $\Waux_{\sf aff}$, $\Jaux_{\sf aff}^{\sim}$
gives $0$. Hence only the converse matters.

After an elementary real 
affine transformation:
\[
u
\,=\,
F(x,y)
\,=\,
x^2
+
{\sf O}_{x,y}(3)
\,\,=\,\,
F_0(y)
+
x\,F_1(y)
+
x^2\,F_2(y)
+
x^3\,F_3(y)
+
x^4\,F_4(y)
+\cdots,
\]
with $F_2(0) = 1$, $F_0(y) = {\sf O}_y(3)$, 
$F_1(y) = {\sf O}_y(2)$. Plug this in~\text{\ding{192}}:
\[
\aligned
0
&
\,\equiv\,
\big(
2\,F_2
+
6x\,F_3
+
{\sf O}_x(2)
\big)\,
\big(
F_{0,yy}
+
x\,F_{1,yy}
+
{\sf O}_x(2)
\big)
-
\big(
F_{1,y}
+
2x\,F_{2,y}
+
{\sf O}_x(2)
\big)^2
\\
&
\,\equiv\,
2\,F_2\,F_{0,yy}
-
\big(F_{1,y}\big)^2
+
x\,
\big[
2\,F_2\,F_{1,yy}
+
6\,F_3\,F_{0,yy}
-
4\,F_{1,y}\,F_{2,y}
\big]
+
{\sf O}_x(2).
\endaligned
\]
Use $F_2(0) \neq 0$ to invert and get:
\[
F_{0,yy}
\,=\,
\mathcal{R}
\cdot
F_{1,y},
\ \ \ \ \ \ \ \ \ \ \ \ \ \ \ \ \ \ \ \ \ \ \ \ \ \
F_{1,yy}
\,=\,
\mathcal{R}
\cdot
F_{0,yy}
+
\mathcal{R}
\cdot
F_{1,y}
\,=\,
\mathcal{R}
\cdot
F_{1,y},
\]
where $\mathcal{R} = \mathcal{R}(y)$ denotes unspecified
functions. From $F_{1,y}(0) = 0$ comes $F_{1,yy}(0) = 0$
and an iteration:
\[
F_{1,yyy}
\,=\,
\mathcal{R}
\cdot
F_{1,y}
+
\mathcal{R}
\cdot
F_{1,yy}
\,=\,
\mathcal{R}
\cdot
F_{1,y},
\,\,\dots\dots,\,\,
F_{1,y^k}
\,=\,
\mathcal{R}
\cdot
F_{1,y},
\,\,\dots\dots,\,\,
\]
yields $F_1(y) \equiv 0$, so $F_{0,yy} \equiv 0$,
whence $F_0(y) \equiv 0$ too. So:
\[
u
\,=\,
x^2
+
\alpha\,x^3
+
\beta\,x^2y
+
{\sf O}_{x,y}(4)
\,=\,
x^2
+
x^2\,
\big(
\underbrace{
\alpha\,x
+
\beta\,y}_{{\sf new}\,\,y}
\big)
+
{\sf O}_{x,y}(4),
\]
since from $2$-nondegeneracy $0 \neq 2 \cdot 2\beta - 0 \cdot
6\alpha$. So:
\[
u
\,=\,
x^2
+
x^2\,y
+
A\,x^4
+
B\,x^3\,y
+
C\,x^2y^2
+
{\sf O}_{x,y}(5).
\]
Then~{\ding{192}}:
\[
0
\,\equiv\,
\Big(
2
+
2\,y
+
{\sf O}_{x,y}(2)
\Big)\,
\big(
2C\,x^2
+
{\sf O}_{x,y}(3)
\big)
-
\Big(
2\,x
+
{\sf O}_{x,y}(2)
\Big)^2
\,\,=\,\,
x^2
\big[
4\,C-4
\big]
+
{\sf O}_{x,y}(3)
\]
forces $C = 1$.

\begin{center}
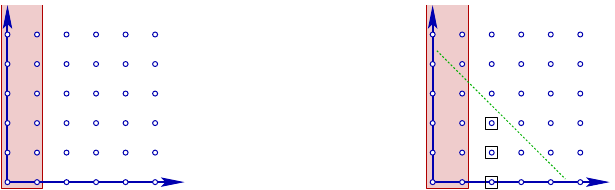
\end{center}

Next, by redefining linearly:
\[
u
\,=\,
x^2
+
x^2\,
\big[
\underbrace{
y+A\,x^2}_{y+A\,u\,=:\,y_\prime}
\big]
+
B\,x^3y
+
x^2y^2
+
{\sf O}_{x,y}(5)
\,\,=\,\,
x^2
+
x^2\,y_\prime
+
B\,x^3\,y_\prime
+
x^2\,y_\prime^2
+
{\sf O}_{x,y_\prime}(5),
\]
we come to:
\[
F
\,=\,
x^2
+
x^2y
+
B\,x^3y
+
x^2y^2
+
{\sf O}_{x,y}(5).
\]

From~{\ding{193}} at $(x,y) = (0,0)$, we kill $0 = 0 - 0 + 
2^2\, 6\, B - 0$. 

We therefore come, after a finite number of affine reductions,
to a suitable form in which $F_{xxx}(0) = 0 = F_{xxxx}(0)$:
\[
F
\,=\,
x^2
+
x^2y
+
x^2y^2
+
{\sf O}_{x,y}(5).
\]

\begin{center}
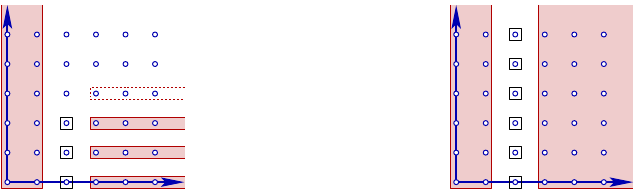
\end{center}

We claim that $F_{x^k} (0) = 0$ for all $k \geqslant 3$.
Indeed, write~{\ding{194}} as $F_{xxxxx} = \mathcal{R}\,
F_{xxx} + \mathcal{R}\, F_{xxxx}$, 
where $\mathcal{R} = \mathcal{R}(x,y)$ is unspecified,
get $F_{xxxxx}(0) = 0$, and
iterate differentiations and substitutions to obtain
$F_{x^k} = \mathcal{R}\, F_{xxx} + 
\mathcal{R}\, F_{xxxx}$ for all $k \geqslant 5$.

We claim that $F_{x^ky}(0) = 0$ for all $k \geqslant 3$.
Indeed, from~{\ding{193}}, solve $F_{xxxy} = \mathcal{R}\, 
F_{xxx} + \mathcal{R}\, F_{xxxx}$, and proceed
similarly.

We claim that $F_{x^ky^\ell}(0) = 0$ for all $k \geqslant 3$ and
$\ell \geqslant 2$. Indeed, from $F_{x^k y^{\ell-1}} = 
\mathcal{R}\, F_{xxx} + \mathcal{R}\, F_{xxxx}$,
differentiate to get:
\[
F_{x^ky^\ell}
\,=\,
\mathcal{R}\,F_{xxx}
+
\mathcal{R}\,F_{xxxy}
+
\mathcal{R}\,F_{xxxx}
+
\mathcal{R}\,F_{xxxxy}
\,\,=\,\,
\mathcal{R}\,F_{xxx}
+
\mathcal{R}\,F_{xxxx}.
\]

So $F(x,y) = x^2\, F_2(y) =: x^2 G(y)$, with
$G(0) = G_y(0) = 1$. Back to~{\ding{192}}
$0 \equiv 2\, G\, x^2\, G_{yy} - \big( 2x\,G_y\big)^2$,
we get:
\[
\aligned
G_{yy}
\,=\,
2!\,\frac{(G_y)^2}{G}
\ \ \ \ \ \ \ \
&
\Longrightarrow
\ \ \ \ \ \ \ \ 
G_{yyy}
\,=\,
2!\,
\frac{2\,G_y\,G_{yy}}{G}
-
2!\,\frac{(G_y)^2\,G_y}{G^2}
\,=\,
3!\,
\frac{(G_y)^3}{G^2}
\\
&
\Longrightarrow
\ \ \ \ \ \ \ \ 
G_{y^k}
\,=\,
k!\,
\frac{(G_y)^k}{G^{k-1}},
\endaligned
\]
whence $G(y) = 1 + y + y^2 + \cdots + y^k + \cdots$ 
and finally after having performed only affine transformations:
\[
u
\,=\,
\frac{x^2}{1-y}.
\qedhere
\]
\endproof

\Section{\bf Affine Invariants via Graph Transforms}
\label{affine-invariants-graph-transforms}
\HEAD{{\ref{affine-invariants-graph-transforms}}.~{\sf Affine 
Invariants via Graph Transforms}
}{
Jo\"el {\sc Merker}, 
Paris-Saclay University, Orsay, France}

In this section, we sketch some
elementary considerations about 
(relative) differential invariants under
$\Aff(\R^3)$. After the prepublication of this paper, 
the article~{\cite{Chen-Merker-2020}}
developed these considerations.
We start with the $1$-dimensional case.

In $\R^2 \ni (x,u)$, the real affine transformation group
$\Aff_2(\R) = \GL_2(\R) \ltimes \R^2$ consists of 
coordinate changes:
\[
\aligned
x'
&
\,=\,
a\,x+b\,y+c,
\\
y'
&
\,=\,
p\,x+q\,y+r,
\endaligned
\]
having nonzero determinant $aq-bp \neq 0$. The {\sl fundamental
equation} expresses how 
graphs are transformed:
\[
p\,x+q\,F(x)+r
\,\equiv\,
F'
\Big(
a\,x+b\,F(x)+c
\Big).
\]
The property of not being a straight line,
namely $F_{xx} \not\equiv 0$, is invariant:
\[
F_{x'x'}'
\,=\,
\frac{(aq-bp)}{\big(a+b\,F_x\big)^3}\,\,
F_{xx}.
\]
Assuming therefore that $F_{xx} \neq 0$ is nowhere vanishing,
whence $F_{x'x'}' \neq 0$ as well,
the Halphen and the Monge 
invariants~{\cite{Halphen-1878, Chen-Merker-2020}}
are well known.

\begin{Theorem}
The Halphen relative invariant whose vanishing characterizes
affine equivalence to $\big\{ u' = x'x' \big\}$ enjoys:
\[
3\,F_{x'x'}'\,F_{x'x'x'x'}'
-
5\,\big(F_{x'x'x'}'\big)^2
\,\,=\,\,
\frac{\big(aq-bp\big)^2}{\big(a+b\,F_x\big)^8}\,
\Big[
3\,F_{xx}\,F_{xxxx}
-
5\,\big(F_{xxx}\big)^2
\Big],
\]
while the Monge relative invariant characterizing the fact
that $\big\{ u = F(x) \big\}$ is contained in a 
nondegenerate conic of
$\R^2$ transforms as:
\reqnomode\usetagform{EngelLie}
\begin{align}
{}
&
9\,\big(F_{x'x'}'\big)^2\,
F_{x'x'x'x'x'}'
-
45\,F_{x'x'}'\,F_{x'x'x'}'\,F_{x'x'x'x'}'
+
40\,\big(F_{x'x'x'}'\big)^3
\,\,=\,\,
\notag
\\
&
\,\,=\,\,
\frac{\big(aq-bp\big)^3}{\big(a+b\,F_x\big)^{12}}\,
\Big[
9\,\big(F_{xx}\big)^2\,
F_{xxxxx}
-
45\,F_{xx}\,F_{xxx}\,F_{xxxx}
+
40\,\big(F_{xxx}\big)^3
\Big].
\tag{\qed}
\end{align}
\end{Theorem}

Next, we pass to the $2$-dimensional case.  As promised, we now
explain how $\Waux_{\sf aff}$ and $\Jaux_{\sf aff}$ can be seen
directly to be affine invariants.  We will even sketch the affine
counterparts of the Levi form, of its kernel field $\mathcal{K}$, of
the nonvanishing function $\laux$, of the slant function $\kaux$, and
of the third-order invariant $\Saux = \overline{\mathcal{L}}_1
(\kaux)$.

In $\R^3 \ni (x,y,u)$, the real affine transformation group
$\Aff_3(\R) = \GLsmall_3 (\R) \ltimes \R^3$ consists of changes of
coordinates:
\[
\aligned
x'
&
\,=\,
a\,x+b\,y+c\,u+d,
\\
y'
&
\,=\,
k\,x+l\,y+m\,u+n,
\\
u'
&
\,=\,
p\,x+q\,y+r\,u+s,
\endaligned
\]
having nonzero Jacobian determinant:
\[
\delta
\,:=\,
\left\vert\!
\begin{array}{ccc}
a & b & c
\\
k & l & m
\\
p & q & r
\end{array}
\!\right\vert
\,\neq\,
0.
\]

We will assume throughout
that such matrices are close to the identity:
\[
\left(\!
\begin{array}{ccc}
a & b & c
\\
k & l & m
\\
p & q & r
\end{array}
\!\right)
\,\,\sim\,\,
\left(\!
\begin{array}{ccc}
1 & 0 & 0
\\
0 & 1 & 0
\\
0 & 0 & 1
\end{array}
\!\right),
\]
so that graphed surfaces $S = \big\{ u = F(x,y) \big\}$
are transformed into similar graphed
surfaces $S' = \big\{ u' = F'(x',y') \big\}$.
This means that by applying the 
$\mathcal{C}^\omega$ implicit function theorem
to the target graphed equation:
\[
p\,x+q\,y+r\,u+s
\,\,=\,\,
F'
\Big(
a\,x+b\,y+c\,u+d,\,\,\,
k\,x+l\,y+m\,u+n
\Big),
\]
the variable $u$ can be solved to recover the first graphed
equation $\big\{ u = F(x,y) \big\}$, that is to say:
\[
u'
\,=\,
F'\big(x',y'\big)
\ \ \ \ \ \ \ \ \ \ 
\Longleftrightarrow
\ \ \ \ \ \ \ \ \ \ 
u
\,=\,
F(x,y).
\]

After preliminary affine normalization, 
we can even assume that $F = {\sf O}_{x,y}(2)$, hence:
\[
F
\,\,\sim\,\,
0,
\ \ \ \ \ \ \ \ \ \ \ \ \ \ \ 
F_x
\,\,\sim\,\,
0,
\ \ \ \ \ \ \ \ \ \ \ \ \ \ \ 
F_y
\,\,\sim\,\,
0.
\]
Then all functions considered will be converging power series
in the two variables $(x,y)$, centered at the origin
$(0,0)$, namely:
\[
F(x,y)
\,\in\,
\R\{x,y\}
\ \ \ \ \ \ \ \ \ \ \ \ \ \ \ \ \ \
\text{and}
\ \ \ \ \ \ \ \ \ \ \ \ \ \ \ \ \ \
F'\big(x',y'\big)
\,\in\,
\R\{x',y'\}.
\]

The {\sl fundamental identity:}
\leqnomode\usetagform{default}
\begin{align}
\label{F-fundamental-F}
p\,x+q\,y+r\,F(x,y)+s
\,\,\equiv\,\,
F'
\Big(
a\,x+b\,y+c\,F(x,y)+d,\,\,\,
k\,x+l\,y+m\,F(x,y)+n
\Big),
\end{align}
holds identically in $\R\{ x, y\}$.

Differentiate this identity with respect to $x$ and to $y$:
\[
\aligned
p+r\,F_x
&
\,\equiv\,
\big(a+c\,F_x\big)\,
F_{x'}'
+
\big(k+m\,F_x\big)\,
F_{y'}',
\\
q+r\,F_y
&
\,\equiv\,
\big(b+c\,F_y\big)\,
F_{x'}'
+
\big(l+m\,F_y\big)\,
F_{y'}'.
\endaligned
\]
To solve for $F_{x'}'$, $F_{y'}'$, a certain $2 \times 2$
determinant appears which depends on the
1\textsuperscript{st} oder jet $J_{x,y}^1F$:
\[
\Lambda
\,:=\,
\Lambda\big(J_{x,y}^1F\big)
\,:=\,
al-bk
+
(cl-bm)\,F_x
+
(am-ck)\,F_y
\,\,\,\sim\,\,
1,
\]
and which is nowhere vanishing, since its value is close to $1$.

Beyond, by differentiating with respect to $x, x$, to
$x, y$, to $y, y$, one solves $F_{x' x'}'$, $F_{x', y'}'$, 
$F_{y', y'}'$ in terms of $J_{x,y}^2F$, and the same determinant
$\Lambda$ appears, as general formulas show
(\cite{Bluman-Kumei-1989, Merker-2008, Chen-Merker-2020}).
The affine invariancy of the Hessian is well known,
and we state a relation that can be verified by a direct
computation.

\begin{Lemma}
One has:
\[
F_{x'x'}'\,F_{y'y'}'
-
\big(F_{x'y'}')^2
\,\,=\,\,
\frac{\delta^2}{\Lambda^4}\,
\Big(
F_{xx}\,F_{yy}
-
\big(F_{xy}\big)^2
\Big).
\eqno\qed
\]
\end{Lemma}

This identity can be abbreviated as:
\[
F_{x'x'}'\,F_{y'y'}'
-
\big(F_{x'y'}')^2
\,\,=\,\,
\nonzero
\cdot
\Big(
F_{xx}\,F_{yy}
-
\big(F_{xy}\big)^2
\Big),
\]
where the generic term `$\nonzero$' denotes various local functions
which are {\em nowhere vanishing}\,\,---\,\,possibly
after shrinking neighborhoods.
Thus, the Hessian is a {\em relative} 
differential invariant under ${\sf SA}_3 (\R)$.

We will make {\em three main hypotheses}, which are meaningful
locally, and which are {\em invariant} under affine transformations.
The first one is:

\begin{Hypothesis}
{\sl The Hessian is degenerate at every point:}
\[
0
\,\equiv\,
F_{xx}\,F_{yy}
-
F_{xy}\,F_{xy}.
\]
\end{Hypothesis}

Not only the Hessian determinant, but also the Hessian matrix
enjoy beautiful invariant properties. Indeed, abbreviate:
\[
\aligned
A(x,y)
&
\,:=\,
a\,x+b\,y+c\,F(x,y)+d,
\\
B(x,y)
&
\,:=\,
k\,x+l\,y+m\,F(x,y)+n,
\\
C(x,y)
&
\,:=\,
p\,x+q\,y+r\,F(x,y)+s,
\endaligned
\]
and differentiate the fundamental 
identity~({\ref{F-fundamental-F}}) once:
\[
\aligned
C_x
&
\,=\,
A_x\,F_{x'}'
+
B_x\,F_{y'}',
\\
C_y
&
\,=\,
A_y\,F_{x'}'
+
B_y\,F_{y'}',
\endaligned
\]
and twice:
\[
\aligned
C_{xx}
&
\,=\,
A_{xx}\,F_{x'}'
+
B_{xx}\,F_{y'}'
\\
&
\ \ \ \ \
+
A_x^2\,F_{x'x'}'
+
2\,A_xB_x\,F_{x'y'}'
+
B_x^2\,F_{y'y'}',
\\
C_{xy}
&
\,=\,
A_{xy}\,F_{x'}'
+
B_{xy}\,F_{y'}'
\\
&
\ \ \ \ \
+
A_xA_y\,F_{x'x'}'
+
\big(
A_xB_y+A_yB_x
\big)\,
F_{x'y'}'
+
B_xB_y\,
F_{y'y'}',
\\
C_{yy}
&
\,=\,
A_{yy}\,F_{x'}'
+
B_{yy}\,F_{y'}'
\\
&
\ \ \ \ \
+
A_y^2\,F_{x'x'}'
+
2\,A_yB_y\,
F_{x'y'}'
+
B_y^2\,F_{y'y'}'.
\endaligned
\]
Introduce the vector fields tangent to $S$ and to $S'$:
\[
\aligned
L_x
&
\,:=\,
\frac{\partial}{\partial x}
+
F_x\,
\frac{\partial}{\partial u},
\ \ \ \ \ \ \ \ \ \ \ \ \ \ \ \ \ \ \ \ \ \ \ \ \ \
L_{x'}
&
\,:=\,
\frac{\partial}{\partial x'}
+
F_{x'}'\,
\frac{\partial}{\partial u'},
\\
L_y
&
\,:=\,
\frac{\partial}{\partial y}
+
F_y\,
\frac{\partial}{\partial u},
\ \ \ \ \ \ \ \ \ \ \ \ \ \ \ \ \ \ \ \ \ \ \ \ \ \
L_{y'}
&
\,:=\,
\frac{\partial}{\partial y'}
+
F_{y'}'\,
\frac{\partial}{\partial u'},
\endaligned
\]
together with their companions, the {\sl horizontal-affine} fields:
\[
\aligned
H_x
&
\,:=\,
\frac{\partial}{\partial x},
\ \ \ \ \ \ \ \ \ \ \ \ \ \ \ \ \ \ \ \ \ \ \ \ \ \
H_{x'}
&
\,:=\,
\frac{\partial}{\partial x'},
\\
H_y
&
\,:=\,
\frac{\partial}{\partial y},
\ \ \ \ \ \ \ \ \ \ \ \ \ \ \ \ \ \ \ \ \ \ \ \ \ \
H_{y'}
&
\,:=\,
\frac{\partial}{\partial y'}.
\endaligned
\]
Although $H_x$, $H_y$ and $H_{x'}$, $H_{y'}$ are {\em not} 
intrinsically related to the geometry of the surfaces
$S$ and $S'$, they will be useful to show that the Hessian
{\em matrices} enjoy invariant properties.
Two natural differential $1$-forms:
\[
\varrho
\,:=\,
du
-
F_x\,dx
-
F_y\,dy
\ \ \ \ \ \ \ \ \ \ \ \ \ \ \ \ \ \
\text{and}
\ \ \ \ \ \ \ \ \ \ \ \ \ \ \ \ \ \
\varrho'
\,:=\,
du'
-
F_{x'}'\,dx'
-
F_{y'}'\,dy',
\]
represent the tangent spaces:
\[
TS
\,=\,
\big\{
\varrho=0
\big\}
\,=\,
\Vect\,
\big(L_x,\,L_y\big)
\ \ \ \ \ \ 
\text{and}
\ \ \ \ \ \ 
TS'
\,=\,
\big\{
\varrho'=0
\big\}
\,=\,
\Vect\,
\big(L_{x'},\,L_{y'}\big).
\]
One may verify that:
\[
\varrho
\,=\,
\mu'\,\varrho',
\]
in terms of the nowhere vanishing function:
\[
\mu'
\,=\,
r
-
c\,F_{x'}'
-
m\,F_{y'}'
\,\,\,\sim\,\,
1.
\]

The proof of the next proposition is left to the reader,
and the reconstitution of appropriate concepts
also, with the hint
of taking inspiration from 
Section~8 of~{\cite{Merker-Pocchiola-Sabzevari-2013-5-CR-II}},
by realizing that the source Hessian matrix
can be written under the appropriate form:
\[
\left(\!
\begin{array}{cc}
F_{xx} & F_{yx}
\\
F_{xy} & F_{yy}
\end{array}
\!\right)
\,\,=\,\,
\left(\!
\begin{array}{cc}
\varrho\big([H_x,\,L_x]\big) & \varrho\big([H_y,\,L_x]\big)
\\
\varrho\big([H_x,\,L_y]\big) & \varrho\big([H_y,\,L_y]\big)
\end{array}
\!\right),
\]
and similarly in the target space:
\[
\left(\!
\begin{array}{cc}
F_{x'x'}' & F_{y'x'}'
\\
F_{x'y'}' & F_{y'y'}'
\end{array}
\!\right)
\,\,=\,\,
\left(\!
\begin{array}{cc}
\varrho'\big([H_{x'},\,L_{x'}]\big) & \varrho'\big([H_{y'},\,L_{x'}]\big)
\\
\varrho'\big([H_{x'},\,L_{y'}]\big) & \varrho'\big([H_{y'},\,L_{y'}]\big)
\end{array}
\!\right).
\]
An alternative direct proof
would be to verify plainly that the shown
matrix identity holds.

\begin{Proposition}
The Hessian matrices in the source space 
$\R_{x,y,u}^3$ and in the target space $\R_{x',y',u'}^3$ enjoy:
\[
\left(\!
\begin{array}{cc}
F_{xx} & F_{yx}
\\
F_{xy} & F_{yy}
\end{array}
\!\right)
\,\,=\,\,
\mu'\,
\left(\!
\begin{array}{cc}
A_x & B_x
\\
A_y & B_y
\end{array}
\!\right)\,
\left(\!
\begin{array}{cc}
F_{x'x'}' & F_{y'x'}'
\\
F_{x'y'}' & F_{y'y'}
\end{array}
\!\right)\,
\left(\!
\begin{array}{cc}
A_x & B_x
\\
A_y & B_y
\end{array}
\!\right)^{\!{\tt t}}.
\eqno\qed
\]
\end{Proposition}

This demonstrates that not only their
(zero) determinants, but also their
ranks are the same!

The most degenerate case occurs when the Hessian matrix is
identically zero.

\begin{Lemma}
The following two conditions
are equivalent for a graphed $\mathcal{C}^\omega$
surface $S = \big\{ u = F(x,y) \big\}$ in $\R^3$.

\smallskip\noindent{\bf (i)}\,
The Hessian {\em matrix} of the graphing function is identically zero: 
\[
F_{xx}
\,\equiv\,
F_{xy}
\,\equiv\,
F_{yx}
\,\equiv\,
F_{yy}
\,\equiv\,
0.
\]

\smallskip\noindent{\bf (ii)}\,
$S$ is affinely equivalent to the flat plane $\big\{ u' = 0\big\}$,
with identically zero graphing function $F' \equiv 0$.\qed

\end{Lemma}

Let us therefore assume that
the rank of the Hessian matrix is $\geqslant 1$.
After some elementary
affine transformation, we come to our second

\begin{Hypothesis}
At every point $F_{xx} \neq 0$.
\end{Hypothesis}

To confirm the invariancy of such a hypothesis, introduce the
nowhere vanishing quantity:
\[
\Upsilon
\,:=\,
\Upsilon\big(J_{x,y}^2F\big)
\,:=\,
\big(l+m\,F_y\big)\,
F_{xx}
-
\big(k+m\,F_x\big)\,
F_{xy}
\,\,\,\sim\,\,
F_{xx}
\,\neq\,
0.
\]

\begin{Lemma}
One has:
\[
F_{x'x'}'
\,=\,
\frac{\delta\,\Upsilon^2}{\Lambda^3}\,
\frac{1}{F_{xx}}.
\eqno\qed
\]
\end{Lemma}

Next, we yet want to exclude the situation where $S = \big\{ u =
F(x,y) \big\}$ is affinely equivalent to $\big\{ u = x^2 \big\}$, a
product of a parabola in $\R_{x,u}^2$ with $\R_y$, and this can be
done by means of a certain relative affine differential
invariant.

\begin{Lemma}
One has:
\[
\frac{F_{x'x'}'\,F_{x'x'y'}'-F_{x'y'}'\,F_{x'x'x'}'}{
\big(F_{x'x'}'\big)^2}
\,\,=\,\,
\frac{F_{xx}}{\Upsilon}\,
\bigg(
\frac{F_{xx}\,F_{xxy}-F_{xy}\,F_{xxx}}{(F_{xx})^2}
\bigg).
\eqno\qed
\]
\end{Lemma}

Similarly as in~{\cite{Pocchiola-2013, 
Merker-Pocchiola-2018, Foo-Merker-2019}}, let
us abbreviate this invariant as:
\[
\Saux_{\sf aff}
\,:=\,
\frac{F_{xx}\,F_{xxy}-F_{xy}\,F_{xxx}}{(F_{xx})^2}.
\]

\begin{Proposition}
\label{Prp-S-not-cylinder}
The following two conditions
are equivalent for a graphed $\mathcal{C}^\omega$
surface $S = \big\{ u = F(x,y) \big\}$ in $\R^3$ satisfying
$F_{xx} \neq 0$ and $0 \equiv F_{xx}\, F_{yy} - F_{xy}^2$.

\smallskip\noindent{\bf (i)}\,
Its relative invariant $\Saux_{\sf aff}$ vanishes identically:
\[
0
\,\equiv\,
F_{xx}\,F_{xxy}
-
F_{xy}\,F_{xxx}.
\]

\smallskip\noindent{\bf (ii)}\,
$S^2 \subset \R^3$ is affinely equivalent to
a cylinder $C^1 \times \R$, with $C^1 \subset \R^2$ a curve.\qed
\end{Proposition}

The proof being done in~{\cite[Section~22]{Chen-Merker-2020}},
we come to our third and last

\begin{Hypothesis}
At every point $F_{xx}\, F_{xxy} - F_{xy}\, F_{xxx} \neq 0$.
\end{Hypothesis}

We mention that 
thanks to the previous formulas,
this numerator of $\Saux_{\sf aff}$ 
and the one of $\Saux_{\sf aff}'$ 
enjoy the transformation rule:
\[
F_{x'x'}'\,F_{x'x'y'}'
-
F_{x'y'}'\,F_{x'x'x'}'
\,\,=\,\,
\frac{\delta^2\,\Upsilon^3}{\Lambda^6}\,
\frac{1}{F_{xx}^3}\,
\Big(
F_{xx}\,F_{xxy}-
F_{xy}\,F_{xxx}
\Big).
\]

\begin{Proposition}
The affinization $\Waux_{\sf aff}$ of 
Pocchiola's invariant $\Waux$ satisfies under an affine equivalence:
\reqnomode\usetagform{EngelLie}
\begin{align}
{}
&
\big(F_{x'x'}'\big)^2\,F_{x'x'x'y'}'
-
F_{x'x'}'\,F_{x'y'}'\,F_{x'x'x'x'}'
+
2\,F_{x'y'}'\,\big(F_{x'x'x'}'\big)^2
-
2\,F_{x'x'}'\,F_{x'x'x'}'\,F_{x'x'y'}'
\,=\,
\notag
\\
&
\ \ \ \ \
\,=\,
\frac{\delta^3\,\Upsilon^6}{(F_{xx})^6\,\Lambda^{10}}\,\,
\Big(
F_{xx}^2\,F_{xxxy}
-
F_{xx}\,F_{xy}\,F_{xxxx}
+
2\,F_{xy}\,F_{xxx}^2
-
2\,F_{xx}\,F_{xxx}\,F_{xxy}
\Big).
\tag{\qed}
\end{align}
\end{Proposition}

Similarly:
\[
\Jaux_{\sf aff}(F')
\,=\,
\frac{\delta^{\sf a}\,\Upsilon^{\sf b}}{\Lambda^{\sf c}}\,\,
\Jaux_{\sf aff}(F),
\]
where ${\sf a}$, ${\sf b}$, ${\sf c}$ are integers,
which can be determined.

\Section{\bf Open Problems}
\label{open-problems}
\HEAD{{\ref{open-problems}}.~{\sf Open Problems}
}{
Jo\"el {\sc Merker}, 
Paris-Saclay University, Orsay, France}

In~{\cite{Olver-2007}}, Olver used
symbolic differential invariants to find a
suprising result:
for a surface $S^2 \subset \R^3$
which is suitably generic,
a {\em single} invariant
generates the whole algebra of differential invariants
under the special affine group
${\sf SA}(\R^3)$,
the {\sl Pick invariant}.
This is also true true under the projective,
conformal, Euclidean groups,
{\em cf.}~{\cite{Arnaldsson-Valiquette-2020}}
and the references therein.

In Spring 2019, we realized that there exist
connections between CR geometry and 
Olver's theory of moving (co)frames.
We were especially interested in the case
where the Hessian matrix:
\[
\left(\!
\begin{array}{cc}
F_{xx} & F_{xy}
\\
F_{yx} & F_{yy}
\end{array}
\!\right)
\]
has constant rank $1$, in view of the analogy with
the much studied class $\mathfrak{C}_{2,1}$ of hypersurfaces
$M^5 \subset \C^3$ having constant Levi rank $1$.
A few e-mail exchanges~{\cite{Olver-2019-personal-communication}}
convinced us that the theory of
moving (co)frames could be applied to {\sl parabolic} surfaces 
$S^2 \subset \R^3$. And generally, we raised

\begin{Question}
\label{Q-differential-invariants-surfaces}
{\sl Study the structure of the full algebra of differential
invariants of real surfaces $\big\{ u = F(x,y) \big\}$
that are {\em parabolic:}}
\[
F_{xx}
\,\neq\,
0
\,\equiv\,
F_{xx}\,F_{yy}
-
F_{xy}^2.
\]
{\sl Study also differential
invariants of surfaces having Hessian of (constant) rank $2$:}
\[
F_{xx}\,
F_{xxy}
-
F_{xy}\,
F_{xxx}
\,\neq\,
0.
\]
\end{Question}

Since the present paper appeared as {\footnotesize\sf
arxiv.org/abs/1903.00889/},
these parabolic surfaces have been studied under
${\sf SA}(\R^3)$ in~{\cite{Chen-Merker-2020}}.
More recently, As of September 2020,
Arnaldsson-Valiquette solved 
Question~{\ref{Q-differential-invariants-surfaces}}
generally, also when the Hessian has rank $2$.

We must formulate a further

\begin{Problem}
\label{Pbm-leaves-homogeneous-models}
{\sl 
From the knowledge of a precise 
{\em branching tree}\,\,---\,\,{\em cf.}
{\cite[Section~2]{Chen-Merker-2020}}\,\,---\,\,of 
differential invariants, after setting up appropriate
recurrence relations for differential
invariants, can one recover the known 
classification~{\cite{Abdalla-Dillen-Vrancken-1997,
Doubrov-Komrakov-Rabinovich-1996,
Eastwood-Ezhov-1999}}
of affinely homogeneous surfaces $S^2 \subset \R^3$?

Specifically, can one determine in which terminal `leaf' 
does each affinely homogeneous model `land'?

Moreover, for each branch, can one determine
a minimal set of 
generating differential invariants (without any homogeneity
assumption)?}
\end{Problem}

Due to the branching process,
trees of differential 
invariants~{\cite{Merker-Nurowski-2020}}
immediately show that homogeneous models 
in different branches 
are {\em not} 
equivalent, 

For instance, 
a partial study of Problem~{\ref{Pbm-leaves-homogeneous-models}} 
was made in~{\cite[Section~23]{Chen-Merker-2020}},
for parabolic surfaces under ${\sf SA}(\R^3)$.
The outcome was that, 
excepting the straight cone 
which can be graphed as $\big\{ u = \frac{1}{2}\, 
\frac{x^2}{1-y} \big\}$, there are no non-cylindrical
special affinely homogeneous parabolic surfaces
$S^2 \subset \R^3$.

The determination of the invariant differential algebras
of parabolic surfaces $S^2 \subset \R^3$ under
the {\em full} affine group ${\sf A}(\R^3)$,
and the determination of all homogeneous models
by such an approach, were done on 
personal notes, were not published
in~{\cite{Chen-Merker-2020}}, but will appear soon. 

\begin{Problem}
{\sl 
Study branching trees, algebras of differential invariants, and
homogeneous models for $3$-dimensional hypersurfaces $H^3 \subset
\R^4$ under the affine group ${\sf A}(\R^4)$, graphed as $\big\{ u =
F(x,y,z) \big\}$, depending on the possible ranks\,\,---\,\,assumed
constant\,\,---\,\,of:}
\[
\Hessian(F)
\,:=\,
\left\vert\!
\begin{array}{ccc}
F_{xx} & F_{xy} & F_{xz}
\\
F_{yx} & F_{yy} & F_{yz}
\\
F_{zx} & F_{zy} & F_{zz}
\end{array}
\!\right\vert.
\]
\end{Problem}

Similar problems can be formulated under other
(finite-dimensional) classical transformation groups 
of $\C^\NN$ or $\R^\NN$: projective;
conformal; Euclidean; {\em etc.}
Also, higher-codimensional submanifolds can be considered,
and so on. To solve in full detail such questions
is probably an undoable infinite task in general dimension 
$\NN \geqslant 1$, because these objects are too wide.

\smallskip

Lastly, we would like to mention an open question 
about Levi degenerate CR geometry
which is inspired from degenerate Affine geometry.
In~{\cite[Section~22]{Chen-Merker-2020}}, 
incited by an anonymous referee, we added the simple

\begin{Observation}
\label{Obs-W-zero-cone}
If $S^2 = \{u = F(x,y)\}$ is a parabolic surface,
{\em i.e.} $F_{xx} \neq 0 = F_{xx}F_{yy} - F_{xy}^2$,
which is not a cylinder, namely
with $0 \neq F_{xx} F_{xxy} - F_{xy} F_{xxx}$
by Proposition~{\ref{Prp-S-not-cylinder}},
then $S^2$ is affinely equivalent to a cone if and only
if $0 \equiv \Waux_{\sf aff}$, that is:
\[
0
\,\equiv\,
\big(F_{xx}\big)^2\,F_{xxxy}
-
F_{xx}\,F_{xy}\,F_{xxxx}
+
2\,F_{xy}\,\big(F_{xxx}\big)^2
-
2\,F_{xx}\,F_{xxx}\,F_{xxy}.
\eqno\qed
\]
\end{Observation}

Recall~{\cite[3-5]{Do-Carmo-2016}}
that a {\sl cone} is a special ruled and even developable 
($=$ parabolic) surface
whose line of striction degenerates to a point, the {\sl vertex} of
the cone.  Putting the vertex at the origin, any 
(smooth) cone can also be
defined parametrically as:
\[
(0,1)\times\R_+^\ast
\ni
(t,v)
\,\,\longmapsto\,\,
v\,\vec{w}(t)
\,\in\,
\R^3,
\]
where $\vec{w}(t) \in P \subset \R^3$ is some space curve
contained in a plane $P \not\ni 0$ not passing through
the origin, which satisfies
$\vec{w}(t) \neq 0 \neq \vec{w}_t(t)$ for all $t$.

Back to a hypersurface
$M^5 \subset \C^3$ in the class $\mathfrak{C}_{2,1}$,
it is known that the Levi kernel bundle 
$\Re\, K^{1,0} M$, of real rank $2$, is Frobenius-integrable,
and that integral $2$-dimensional manifolds are
{\em holomorphic curves}.
So every such $M^5$ is {\em foliated}
by complex curves, hence can be locally parametrized as:
\[
\Phi\colon\ \ \ \ \
\aligned
\C\times\R\times\R\times\R
&
\,\,\,\longrightarrow\,\,\,
M^5
\subset
\C^3
\\
(\tau,r,s,t)
&
\,\,\,\longmapsto\,\,\,
\big(
z(\tau,r,s,t),\,\,
\zeta(\tau,r,s,t),\,\,
w(\tau,r,s,t)
\big),
\endaligned
\]
with $\Phi$ analytic as $M$ is, and moreover,
{\em holomorphic} with respect to $\tau$,
that is $0 \equiv z_{\overline{\tau}} 
\equiv \zeta_{\overline{\tau}} \equiv w_{\overline{\tau}}$.

\begin{Question}
\label{Q-Waux-zero-holomorphic-cone}
{\sl 
Is there a natural geometric interpretation of the identical
vanishing $0 \equiv \Waux_0$ of Pocchiola's 
(first) relative invariant about the parametrization
$\Phi(\tau, r, s, t)$,
which would be analogous, in some sense,
to Observation~{\ref{Obs-W-zero-cone}}?}
\end{Question}

Recently, Merker-Nurowski~{\cite{Merker-Nurowski-2020}}
considered {\sl para-CR structure} associated 
to CR hypersurfaces $M^5 \in \mathfrak{C}_{2,1}$.
Such para-CR structures can be seen as systems of two PDEs:
\[
z_y
=
G(x,y,z,z_x,z_{xx})
\quad\&\quad 
z_{xxx}
=
H(x,y,z,z_x,z_{xx}),
\quad
\text{for}\,\,
z=z(x,y),
\]
with complete integrability $D_x D_x D_x G = D_y H$, in terms of two
real $\mathcal{C}^\omega$
functions $G = G(x,y,z,p,r)$ and $H = H(x,y,z,p,r)$ satisfying
$G_r\equiv 0$ and $G_{pp} \neq 0$ to insure constant Levi rank $1$ and
$2$-nondegeneracy.

It can be verified that $0 \equiv \Waux_0$ then translates as the
vanishing of one among three primary relative para-CR differential
invariants:
\[
0
\,\equiv\,
2G_{ppp}+G_{pp}H_{rr}.
\]
Under this assumption, it was discovered 
in~{\cite{Merker-Nurowski-2020}}
that such geometric structures {\em embed}
(in a certain sense) into the equivalence
problem of 3\textsuperscript{rd} order ODEs
under {\em contact} transformations. 

So various geometric perspectives\big/interpretations 
can be expected to answer 
Question~{\ref{Q-Waux-zero-holomorphic-cone}}. 

\Section{\bf Smoothness Assumption Improvement}
\label{smoothness-assumption-improvement}
\HEAD{{\ref{smoothness-assumption-improvement}}.~{\sf Smoothness 
Assumption Improvement}
}{
Jo\"el {\sc Merker}, 
Paris-Saclay University, Orsay, France}

In~{\cite[3.1]{Isaev-2011}}, it was shown that
any $\mathcal{C}^\infty$ tube hypersurface
$\big\{ \Re\, w = F(\Re\,z_1, \dots, \Re\, z_n) \big\}$
in $\C^{n+1}$ which is (locally) biholomorphic
to the sphere $\big\{ \Re\, w = (\Re\, z_1)^2 + \cdots + 
(\Re\, z_n)^2 \big\}$\,\,---\,\,in its unbounded 
representation\,\,---\,\,is in fact real analytic, 
$\mathcal{C}^\omega$.
However, this is untrue for general spherical 
hypersurfaces, that are neither tube nor even rigid,
{\em see}~{\cite[3.3]{Isaev-2011}}.
Let us formulate and prove an analogous

\begin{Proposition}
Every $\mathfrak{C}_{2,1}$ tube hypersurface
$\big\{ \Re\, w = F(\Re\, z, 
\Re\, \zeta) \big\}$
of class $\mathcal{C}^5$, {\em e.g.} of class
$\mathcal{C}^\infty$, 
which is biholomorphic to the model
$\big\{ \Re\, w = \frac{(\Resmall\, z)^2}{
1-\Resmall\, \zeta} \big\}$,
is in fact $\mathcal{C}^\omega$.
\end{Proposition}

Hence Theorem~{\ref{Thm-main}}, stated in the $\mathcal{C}^\omega$
category, holds in fact for $M$ of class $\mathcal{C}^5$.

\begin{center}
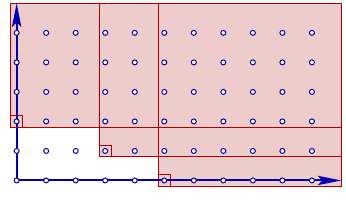
\end{center}

\proof
Denote $x := \Re\, z$, $y := \Re\, \zeta$, $u := \Re\, w$.
According to Theorem~{\ref{Theorem-LC-J-W-zero}},
the graphing function $F$ satisfies $3$ PDEs.
Since $F_{xx} \neq 0$, we can rewrite these
PDEs in solved form as:
\[
\aligned
F_{yy}
&
\,=\,
\frac{F_{xy}^2}{F_{xx}},
\\
F_{xxxy}
&
\,=\,
\frac{F_{xy}}{F_{xx}}\,
F_{xxxxx}
+
2\,\frac{F_{xxx}\,F_{xxy}}{F_{xx}}
-
2\,\frac{F_{xy}\,F_{xxx}^2}{F_{xx}^2},
\\
F_{xxxxx}
&
\,=\,
\frac{45}{9}\,
\frac{F_{xxx}\,F_{xxxx}}{F_{xx}}
-
\frac{40}{9}\,\frac{F_{xxx}^3}{F_{xx}^2}.
\endaligned
\]
Although the computational task is not straightforward,
one can verify that this PDE system is compatible.
In fact, for every $j + k \leqslant 5$, with:
\[
(j,k)
\,\neq\,
\aligned
{}
&
(0,1),\ \ \ (1,1),\ \ \ (2,1),
\\
{}
&
(0,0),\ \ \ (1,0),\ \ \ (2,0),\ \ \ (3,0),\ \ \ (4,0),
\endaligned
\]
one can compute certain {\em uniquely defined} 
right-hand sides:
\[
F_{x^jy^k}
\,=\,
\mathcal{R}_{j,k}
\left(
\aligned
{}
&
F_{xy},\ \ \
F_{xxy},
\\
{}
&
\ \ \ \ \ \ \ \ \ \ \ \ 
F_{xx},\ \ \ F_{xxx},\ \ \ F_{xxxx}
\endaligned
\right).
\]

Equivalently, in the jet space equipped with 
coordinates $u_{x^j y^k}$, 
one can verify that the two total
differentiation operators, restricted to the 
PDE system, namely:
\[
\aligned
&
\ \ \ \ \ \ \ \ \ \ \ \ \ \ 
+
u_{xy}\,\frac{\partial}{\partial u_y}
+
u_{xxy}\,\frac{\partial}{\partial u_{xy}}
+
\mathcal{R}_{3,1}\,\frac{\partial}{\partial u_{xxy}}
\\
D_x
&
\,:=\,
\frac{\partial}{\partial x}
+
u_x\,\frac{\partial}{\partial u}
+
u_{xx}\,\frac{\partial}{\partial u_x}
+
u_{xxx}\,\frac{\partial}{\partial u_{xx}}
+
u_{xxxx}\,\frac{\partial}{\partial u_{xxx}}
+
\mathcal{R}_{5,0}\,\frac{\partial}{\partial u_{xxxx}},
\endaligned
\]
and:
\[
\aligned
&
\ \ \ \ \ \ \ \ \ \ \ \ \ \ 
+
\mathcal{R}_{0,2}\,\frac{\partial}{\partial u_y}
+
\mathcal{R}_{1,2}\,\frac{\partial}{\partial u_{xy}}
+
\mathcal{R}_{2,2}\,\frac{\partial}{\partial u_{xxy}},
\\
D_y
&
\,:=\,
\frac{\partial}{\partial y}
+
u_y\,\frac{\partial}{\partial u}
+
u_{yx}\,\frac{\partial}{\partial u_x}
+
u_{yxx}\,\frac{\partial}{\partial u_{xx}}
+
\mathcal{R}_{3,1}\,\frac{\partial}{\partial u_{xxx}}
+
\mathcal{R}_{4,1}\,\frac{\partial}{\partial u_{xxxx}},
\endaligned
\]
{\em commute one with another}, {\em cf.}~{\cite{Merker-2008}}
in a general context.

Consequently, similarly as for classical second order 
ODEs, in the jet submanifold equipped with coordinates:
\[
\left(
\aligned
{} 
&
\ \ \ \ \ \ \ \ \
u_y,\,\,
u_{xy},\,\,
u_{xxy},
\\
{}
&
x,y,\,\,
u,\,\,
u_x,\,\,
u_{xx},\,\,
u_{xxx},\,\,
u_{xxxx}
\endaligned
\right),
\]
the Frobenius theorem applies, and provides a general solution:
\[
F
\,=\,
F(x,y)
\,=\,
\mathcal{F}
\left(
\aligned
{} 
&
\ \ \ \ \ \ \ \ \
u_{0,1}^0,\,\,
u_{1,1}^0,\,\,
u_{2,1}^0,
\\
{}
&
x,y,\,\,
u_{0,0}^0,\,\,
u_{1,0}^0,\,\,
u_{2,0}^0,\,\,
u_{3,0}^0,\,\,
u_{4,0}^0
\endaligned
\right),
\]
which is {\em analytic}, {\em i.e.} $\mathcal{C}^\omega$,
{\em because the $\mathcal{R}_{j,k}$ are $\mathcal{C}^\omega$},
with respect to {\em all} arguments, 
including the
parameters $u_{j,k}^0$, with the initial
conditions property that:
\[
\footnotesize
\aligned
F_y(0,0)
&
=
u_{0,1}^0,
\ \ \ \ \ \ \ 
F_{xy}(0,0)
=
u_{1,1}^0,
\ \ \ \ \ \ \ 
F_{xxy}(0,0)
=
u_{2,1}^0,
\\
F(0,0)
&
=
u_{0,0}^0,
\ \ \ \ \ \ \ 
F_x(0,0)
=
u_{1,0}^0,
\ \ \ \ \ \ \ 
F_{xx}(0,0)
=
u_{2,0}^0,
\ \ \ \ \ \ \ 
F_{xxx}(0,0)
=
u_{3,0}^0,
\ \ \ \ \ \ \ 
F_{xxxx}(0,0)
=
u_{4,0}^0,
\endaligned
\]

In conclusion, if $F \in \mathcal{C}^5$ satisfies
the $3$ PDEs in question,
its representation as $F = \mathcal{F}$ with such
constants
$u_{j,k}^0 = F_{x^jy^k}(0,0)$ shows that 
$F \in \mathcal{C}^\omega$.
\endproof






\vfill\end{document}

%% file: macros.tex
\newtheorem{Theorem}[equation]{Theorem}

\newtheorem{Proposition}[equation]{Proposition}

\newtheorem{Lemma}[equation]{Lemma}

\newtheorem{Corollary}[equation]{Corollary}

\newtheorem{Observation}[equation]{Observation}

\theoremstyle{definition}

\newtheorem{Hypothesis}[equation]{Hypothesis}

\newtheorem{Definition-Notation}[equation]{D\'efinition-Notation}
\newtheorem{Notation}[equation]{Notation}

\newtheorem{Question}[equation]{Question}

\newtheorem{Problem}[equation]{Problem}

\newcommand{\C}{\mathbb{C}}

\newcommand{\K}{\mathbb{K}}

\renewcommand{\P}{\mathbb{P}}

\newcommand{\R}{\mathbb{R}}

\newcommand{\NN}{\text{\sc n}}

\newcommand{\kaux}{{\text{\usefont{T1}{qcs}{m}{sl}k}}}
\newcommand{\laux}{{\text{\usefont{T1}{qcs}{m}{sl}l}}}

\newcommand{\xaux}{{\text{\usefont{T1}{qcs}{m}{sl}x}}}
\newcommand{\yaux}{{\text{\usefont{T1}{qcs}{m}{sl}y}}}
\newcommand{\zaux}{{\text{\usefont{T1}{qcs}{m}{sl}z}}}

\newcommand{\Aaux}{{\text{\usefont{T1}{qcs}{m}{sl}A}}}

\newcommand{\Faux}{{\text{\usefont{T1}{qcs}{m}{sl}F}}}

\newcommand{\Haux}{{\text{\usefont{T1}{qcs}{m}{sl}H}}}
\newcommand{\Iaux}{{\text{\usefont{T1}{qcs}{m}{sl}I}}}
\newcommand{\Jaux}{{\text{\usefont{T1}{qcs}{m}{sl}J}}}

\newcommand{\Paux}{{\text{\usefont{T1}{qcs}{m}{sl}P}}}

\newcommand{\Raux}{{\text{\usefont{T1}{qcs}{m}{sl}R}}}
\newcommand{\Saux}{{\text{\usefont{T1}{qcs}{m}{sl}S}}}

\newcommand{\Waux}{{\text{\usefont{T1}{qcs}{m}{sl}W}}}

\definecolor{blue}{cmyk}{1.,1.,0.,0.63}
\definecolor{red}{cmyk}{0.,1.,1.,0.63}
\definecolor{green}{cmyk}{1.,0.,1.,0.63}
\definecolor{black}{cmyk}{1.,1.,1.,1.}

\newcommand{\blue}{\textcolor{blue}}

\newcommand{\red}{\textcolor{red}}

\makeatletter
\renewcommand{\@fnsymbol}[1]
{\ensuremath{\ifcase#1\or $*$\or $**$\or $***$\or $****$\or $*****$
\else\@ctrerr\fi}}
\makeatother

\newcommand{\HEAD}[2]{%
\pagestyle{fancy}
\fancyhead[RO]{\tiny\sf\thepage}
\fancyhead[CO]{{\tiny\sf #1}}
\fancyhead[LE]{\tiny\sf\thepage}
\fancyhead[CE]{{\tiny\sf #2}}
\fancyfoot{}}

\numberwithin{equation}{section}

\newcommand{\Section}[1]{
\renewcommand{\thesection}{\bf\arabic{section}}
\section{#1}
\renewcommand{\thesection}{\arabic{section}}}

\newcommand{\style}[1]{\text{\footnotesize{\sf #1}}}

\newcommand{\stylesmall}[1]{{\sf #1}}

\newcommand{\Aff}{\style{Aff}}
\newcommand{\Affsmall}{\stylesmall{Aff}}

\renewcommand{\arcsin}{\style{arcsin}}

\newcommand{\arcsinh}{\style{arcsinh}}

\newcommand{\Bihol}{\style{Bihol}}
\newcommand{\Biholsmall}{\stylesmall{Bihol}}

\newcommand{\Cartansmall}{\stylesmall{Cartan}}

\renewcommand{\dim}{\style{dim}}

\newcommand{\GL}{\style{GL}}
\newcommand{\GLsmall}{\stylesmall{GL}}

\newcommand{\Halphensmall}{\stylesmall{Halphen}}

\newcommand{\Hessian}{\style{Hessian}}

\renewcommand{\Im}{\style{Im}}

\newcommand{\Levi}{\style{Levi}}

\renewcommand{\lim}{\style{lim}}

\renewcommand{\log}{\style{log}}

\renewcommand{\mod}{\style{mod}}

\newcommand{\nonzero}{\style{nonzero}}

\renewcommand{\Re}{\style{Re}}
\newcommand{\Resmall}{\stylesmall{Re}}

\newcommand{\SO}{\style{SO}}

\newcommand{\Sym}{\style{Sym}}

\newcommand{\tubesmall}{\stylesmall{tube}}

\newcommand{\Vect}{\style{Vect}}

\newcommand{\Hall}{\Hall}

\newcommand{\vf}{\vfill\end{document}}

\makeatletter
\def\hlinethick#1{%
\noalign{\ifnum0=`}\fi\hrule \@height #1 \futurelet
\reserved@a\@xhline}
\makeatother

\newlength{\myskip}


%% file: 1-F-x-xi-destroy.pdf_t
\begin{picture}(0,0)%
\includegraphics{1-F-x-xi-destroy.pdf}%
\end{picture}%
\setlength{\unitlength}{4144sp}%
\begingroup\makeatletter\ifx\SetFigFont\undefined%
\gdef\SetFigFont#1#2#3#4#5{%
  \reset@font\fontsize{#1}{#2pt}%
  \fontfamily{#3}\fontseries{#4}\fontshape{#5}%
  \selectfont}%
\fi\endgroup%
\begin{picture}(4669,1444)(844,-1468)
\put(2164,-1379){\makebox(0,0)[lb]{\smash{{\SetFigFont{8}{9.6}{\familydefault}{\mddefault}{\updefault}{\color[rgb]{0,0,0}$x$}%
}}}}
\put(936,-147){\makebox(0,0)[lb]{\smash{{\SetFigFont{8}{9.6}{\familydefault}{\mddefault}{\updefault}{\color[rgb]{0,0,0}$\xi$}%
}}}}
\put(5404,-1379){\makebox(0,0)[lb]{\smash{{\SetFigFont{8}{9.6}{\familydefault}{\mddefault}{\updefault}{\color[rgb]{0,0,0}$x$}%
}}}}
\put(4176,-147){\makebox(0,0)[lb]{\smash{{\SetFigFont{8}{9.6}{\familydefault}{\mddefault}{\updefault}{\color[rgb]{0,0,0}$\xi$}%
}}}}
\end{picture}%

%% file: 2-F-x-xi-destroy.pdf_t
\begin{picture}(0,0)%
\includegraphics{2-F-x-xi-destroy.pdf}%
\end{picture}%
\setlength{\unitlength}{4144sp}%
\begingroup\makeatletter\ifx\SetFigFont\undefined%
\gdef\SetFigFont#1#2#3#4#5{%
  \reset@font\fontsize{#1}{#2pt}%
  \fontfamily{#3}\fontseries{#4}\fontshape{#5}%
  \selectfont}%
\fi\endgroup%
\begin{picture}(4849,1444)(664,-1468)
\put(5404,-1379){\makebox(0,0)[lb]{\smash{{\SetFigFont{8}{9.6}{\familydefault}{\mddefault}{\updefault}{\color[rgb]{0,0,0}$x$}%
}}}}
\put(756,-147){\makebox(0,0)[lb]{\smash{{\SetFigFont{8}{9.6}{\familydefault}{\mddefault}{\updefault}{\color[rgb]{0,0,0}$\xi$}%
}}}}
\put(4176,-147){\makebox(0,0)[lb]{\smash{{\SetFigFont{8}{9.6}{\familydefault}{\mddefault}{\updefault}{\color[rgb]{0,0,0}$\xi$}%
}}}}
\end{picture}%

%% file: 0-2-3-1-5-0.pdf_t
\begin{picture}(0,0)%
\includegraphics{0-2-3-1-5-0.pdf}%
\end{picture}%
\setlength{\unitlength}{4144sp}%
\begingroup\makeatletter\ifx\SetFigFont\undefined%
\gdef\SetFigFont#1#2#3#4#5{%
  \reset@font\fontsize{#1}{#2pt}%
  \fontfamily{#3}\fontseries{#4}\fontshape{#5}%
  \selectfont}%
\fi\endgroup%
\begin{picture}(2624,1552)(774,-1591)
\put(1895,-1549){\makebox(0,0)[lb]{\smash{{\SetFigFont{6}{7.2}{\familydefault}{\mddefault}{\updefault}{\color[rgb]{0,0,0}\red{$(5,\!0)$}}%
}}}}
\put(789,-1099){\makebox(0,0)[lb]{\smash{{\SetFigFont{6}{7.2}{\familydefault}{\mddefault}{\updefault}{\color[rgb]{0,0,0}\red{$(0,\!2)$}}%
}}}}
\put(1458,-1320){\makebox(0,0)[lb]{\smash{{\SetFigFont{6}{7.2}{\familydefault}{\mddefault}{\updefault}{\color[rgb]{0,0,0}\red{$(3,\!1)$}}%
}}}}
\put(3276,-1358){\makebox(0,0)[lb]{\smash{{\SetFigFont{8}{9.6}{\familydefault}{\mddefault}{\updefault}{\color[rgb]{0,0,0}\blue{$j$}}%
}}}}
\end{picture}%